 \documentclass[10pt]{article}   

\usepackage{amsthm,amsmath,amssymb,mathrsfs,cite,amscd,curvesls}
\usepackage[dvips]{graphicx}
\numberwithin{equation}{section}

\newcommand{\p}{\partial}

\newcommand{\half}{\tfrac12}

\newcommand{\CO}{\mathcal{O}}

\renewcommand{\geq}{\geqslant}
\renewcommand{\leq}{\leqslant}

\newcommand{\N}{\mathbb{N}}

\newcommand{\R}{\mathbb{R}}
\newcommand{\Z}{\mathbb{Z}}

\newcommand{\li}{\textrm{Li}_2}
\newcommand{\Li}{{\rm Li}}

 \newcommand{\bs}{\begin{split}}
 \newcommand{\be}{\begin{equation}}
 \newcommand{\es}{\end{split}} 
 \newcommand{\ee}{\end{equation}}

\theoremstyle{plain}

\newtheorem{lem}[equation]{Lemma}

\newtheorem{thm}[equation]{Theorem}
\newtheorem{conj}[equation]{Conjecture}

\theoremstyle{remark}

\begin{document}

\begin{titlepage} 
  \linespread{1.8}
  \title{\Large \bf On asymptotics, Stirling numbers, Gamma function
  and polylogs}
  \author{Daniel B. Gr\"unberg \\ {\em \small MPI Bonn, Vivatsgasse 7,
  53111 Bonn, Germany} \\ \small grunberg@mccme.ru} 
  \date{Dec 2005}
  \maketitle
  \begin{abstract}   \large
    We apply the Euler--Maclaurin formula to find the asymptotic
    expansion of the sums $\sum_{k=1}^n (\log k)^p / k^q$, ~$\sum k^q
    (\log k)^p$, ~$\sum (\log k)^p /(n-k)^q$, ~$\sum 1/k^q (\log k)^p
    $ in closed form to arbitrary order ($p,q \in\N$). The expressions
    often simplify considerably and the coefficients are recognizable
    constants.  The constant terms of the asymptotics are either
    $\zeta^{(p)}(\pm q)$ (first two sums), 0 (third sum) or yield
    novel mathematical constants (fourth sum). This allows numerical
    computation of $\zeta^{(p)}(\pm q)$ faster than any current
    software.  One of the constants also appears in the expansion of
    the function $\sum_{n\geq 2} (n\log n)^{-s}$ around the
    singularity at $s=1$; this requires the asymptotics of the
    incomplete gamma function.  The manipulations involve polylogs for
    which we find a representation in terms of Nielsen integrals, as
    well as mysterious conjectures for Bernoulli numbers. Applications
    include the determination of the asymptotic growth of the Taylor
    coefficients of $(-z/\log(1-z))^k$.  We also give the asymptotics
    of Stirling numbers of first kind and their formula in terms of
    harmonic numbers.
  \end{abstract}
\vspace{1cm}
{\small To appear in: Results in Mathematics.}
  \thispagestyle{empty}
\end{titlepage}

\section{Introduction}

This paper is about concrete mathematics. It gathers several results
about asymptotic theory, half of which are obtained from the
Euler--Maclaurin formula.  A few by-products offer themselves, such as
the asymptotics of the incomplete gamma function, or the study of the
complex function $\sum_{n\geq 2} (n\log n)^{-s}$ with a singularity at
$s=1$, or representations of polylogs in terms of Nielsen integrals,
or properties of Stirling numbers, or some identities about Bernoulli
numbers.  We also summarise three ways of obtaining the asymptotic
growth of the Taylor coefficients of $(-z/\log(1-z))^k$.  Much of the
contents may not be new -- let alone ground-breaking, but the interest
of the paper lies in the way all these objects tie the knot and pop up
by studying a few simple problems; it will offer some surprises to the
curious and hands-on mathematician.

To begin with, we recall the Euler--Maclaurin formula:
$$
\sum_{k=a}^{n-1} f(k) = \int_a^n f(x) dx -\half [f(n) -f(a)]+
\sum_{k=1}^m \frac{B_{2k}}{(2k)!} [f^{(2k-1)}(n) -f^{(2k-1)}(a)]
+\textrm{error}
$$
where the error term is $O(\frac{1}{(2\pi)^{2m}})\int_a^n
|f^{(2m)}(x)| dx $.  The values of $\frac{B_{2k}}{(2k)!}$ are
$\frac{1}{12}, -\frac{1}{720}, \frac{1}{30240},\dots$.  

We shall be interested in the limit of large $n$, keeping $a$ fixed.
When ordering the terms in decreasing orders of $n$, the quantity $
-\sum_{k=1}^m \frac{B_{2k}}{(2k)!} f^{(2k-1)}(a)$ will contribute to
the constant term.  The constant term will be exact when all orders
have been taken into account (ie. $m\to\infty$).  Since this means
adding always bigger chunks ($B_{2k} \sim (2k)^{2k}$), we would end up
with an infinite value for the constant term.  In practice, the exact
value of the constant term has to be computed from another approach.
However, the formal infinite sum involving Bernoulli numbers appears
most useful, as it behaves linearly: adding two such sums (from the
asymptotics of $H_n$ and $H_n^{(2)}$, say) will stand for a constant
whose exact value is the sum of the two exact values of the respective
constant terms.  

We shall use this trick in \underline{section 2} to write down the
exact constants hiding behind formal sums.  They will prove useful in
subsequent sections to derive the coefficients in the asymptotic
expansions (for large $n$) of the four sums that we consider in
sections 4,5,6,7 respectively:
$$
\sum_{k=1}^n \frac{(\log k)^p}{ k^q}, \qquad \sum_{k=1}^n k^q (\log
k)^p, \qquad \sum_{k=1}^{n-1} \frac{(\log k)^p}{(n-k)^q}, \qquad
\sum_{k=2}^{n-1} \frac{1}{k^q (\log k)^p} 
$$
for $p,q \in\N$.  We shall write their asymptotics in closed form
to arbitrary order of $n$.  In particular, we can write down
$\zeta^{(p)}(\pm q)$ and the Stieltjes constants $\gamma_p$ as formal
sums over rational numbers.  In this formal sense,  $\gamma_p =(-1)^p
\zeta^{(p)}(1)$.

The coefficients in the asymptotic expansions often contain Stirling
numbers of the first kind, or their close relative which we denote by
$S_{r,s,t} := \sum \frac{1}{i_1\cdots i_r}$ (sum over all integers
$i_j$ such that $s \leq i_1<...< i_r \leq t$).  \underline{Section~3}
expresses these numbers in terms of harmonic numbers, which allows a
rapid deduction of their asymptotics to arbitrary order.  The formula
can be inverted to express harmonic numbers in terms of Stirling
numbers.

The asymptotic expansion of the four sums, presented in
\underline{sections 4,5,6,7}, can be easily derived from the
Euler--Maclaurin formula for the first two sums but involves intricate
algebra for the latter two.  In those cases, the expansion was first
found empirically using the asymp$_k$ trick (appendix).  The
coefficients are all rational numbers except for the constant terms
($\zeta^{(p)}(\pm q)$ for the first two sums, unknown constants for
the fourth sum).  For the third sum, the constant term vanishes but
$\zeta^{(p)}(-q)$ occurs at higher orders (irrational).  The asymp$_k$
trick gives us sufficient digits of a coefficient $c$; we then can use
the PARI software to find a vanishing integer linear combination of
$1,c,\zeta'(7)$, say, if one suspected there was a $\zeta'(7)$ hiding
behind $c$.  The proper linear combination requires often guesswork.
As an application, knowing a large number of terms of the asymptotic
expansion of the first two sums allows one to compute the constant
term $\zeta^{(p)}(\pm q)$ to arbitrary precision more rapidly than any
current mathematical software; the asymp$_k$ trick can also enhance
speed.

As an application, we derive in \underline{section 8} the asymptotic
growth of the coefficients in the Taylor expansion of
$(-z/\log(1-z))^k$, via a convolution from the ansatz at $k=1$ (the
latter known to P\'olya).  The result appeared in two other contexts
(\cite{N-61} and \cite{FO-90}) which we recapitulate for the
interested reader.

As advertised, the fourth sum (section 7) gives birth to a 2d-array of
unknown mathematical constants, $C_{p,q}$, that converge to the values
of $1/(2^q (\log 2)^p)$ when $p,q\to\infty$; only $C_{1,0}$ and
$C_{1,1}$ have appeared (indirectly) before in the literature.
\underline{Section 9} verifies that $C_{1,1}$, which occurs in
$\sum_{k=1}^n \frac{1}{k \log k} \approx \log\log n +C_{1,1} +O(\tfrac
1{n \log n})$, also occurs in the constant term of the asymptotic
expansion of the following complex function around its singularity at $s=1$:
$$
 \sum_{n=2}^\infty \frac{1}{(n\log n)^s} \approx -\log(s-1) +C_{1,1}-\gamma
 +O\big((s-1)\log^2(s-1)\big) \qquad \textrm{as }~ s\to 1.
$$
This involves the asymptotics of the incomplete gamma function.

In order to prove the asymptotic expansion of the third sum (section
6) via the Euler--Maclaurin formula, one needs to track down
surprising cancellations. The manipulations involve a particular
representation of polylogs by Nielsen integrals
$\mathfrak{S}_{1,p}(x)$ presented in \underline{section~10}:
$$
\Li_j(1-x) = \sum_{r=0}^{j-1} \Big( \zeta(j-r) -
\mathfrak{S}_{1,j-r-1}(x) \Big) \frac{\log^{r}(1-x)}{r!},
$$
wherein the term with $\zeta(1)$ should be dropped.  The
generalised polylogs $\Li_{s_1,\dots,s_k}(x) := \sum
\frac{x^{n_1}}{n_1^{s_1} \dots n_k^{s_k}}$ (sum over integers $n_j$
with $n_1>\dots >n_k>0 $) give rise to the Nielsen integrals:
$\mathfrak{S}_{k,p}(x)= \Li_{k+1,1^{p-1}}(x)$ (the subscript $1^{p-1}$
stands for $p-1$ times 1). Thus, the representation can be rewritten
as
$$
\Li_j(1-x) = \sum_{r=0}^{j-1} \Big( \Li_{j-r}(1) -
\Li_{2,1^{j-r-2}}(x) \Big) (-1)^r ~\Li_{1^r}(x),
$$
wherein the term with $\Li_1(1)$ should be dropped.  
Proving the asymptotics of the third sum for $p=2$ entails two curious
representations of Nielsen integrals, (\ref{eq:S11}) and
(\ref{eq:S12}), which themselves boil down to the following bizarre
identities for Bernoulli numbers: For $n$ a positive integer, $n\geq
2$,
\begin{align*}
\sum_{r=1}^{n-1} \frac{(-1)^r B_r}{r} \sum_{l=r}^n \frac{(-1)^l}{l}
\binom{n-1}{n-l} &= -\frac{1}{n^2}   \\
\sum_{r=1}^{n-1} \frac{(-1)^r B_r}{r} \Big( \sum_{l=r}^n (-1)^l
\binom{n}{l} H_{l-1} +\frac{1}{r}+\frac{1}{n-r} \Big) &= H_{n-1}^{(2)}
+\frac{1}{n} H_{n-1}.  
\end{align*}
Proving the asymptotics for higher $p$, one gets a further such
identity, and a whole tower can be built up. The first identity is
easy to prove, but the second has resisted our best efforts (and those
of experts).

\section{Formal sums of Bernoulli numbers and zeta-values }

We start with formal infinite sums involving Bernoulli numbers.  The
notation is formal because the sums diverge ($B_{2k} \sim
(2k)^{2k}$).  Nevertheless, they are useful as one can recognize
constant terms from the expressions $\sum \frac{B_{2k}}{(2k)!} c_k$ in
the Euler--Maclaurin formula.  We shall use such (diverging)
expressions to recognize constants in future applications of the
Euler--Maclaurin formula.

\begin{lem} \label{lemma0}
In formal notation:
  \begin{align}
\gamma &= \half + \sum_1 \frac{B_{2k}}{2k} \\
\zeta(2) &= \tfrac{3}{2} + \sum B_{2k} \\
\zeta(3) &= 1 + \sum B_{2k} \frac{2k+1}{2} \\
\zeta(4) &= \tfrac{5}{6} + \sum B_{2k} \frac{(2k+1)(2k+2)}{2\cdot 3} \\
\zeta(i) &=\frac{i+1}{2i-2} +\sum \frac{B_{2k}}{(2k)!}
\frac{(2k+i-2)!}{(i-1)!} \label{eq:zeta(i)}
  \end{align}
  \begin{align}
\half \log(2\pi)= -\zeta'(0) &= 1-  \sum \frac{B_{2k}}{(2k)(2k-1)} \\
\tfrac{1}{12}-\zeta'(-1) &= \tfrac{1}{4} +\sum_2
\frac{B_{2k}}{(2k)(2k-1)(2k-2)} \\ 
-\zeta'(-2) &= \tfrac{1}{36} -2!\sum_2 \frac{B_{2k}}{(2k)\dots(2k-3)} \\
-\tfrac{11}{720}-\zeta'(-3) &= -\tfrac{1}{48} +3!\sum_3
\frac{B_{2k}}{(2k)\dots(2k-4)} \\ 
\tfrac{B_{q+1}}{q+1} H_q -\zeta'(-q) &= \tfrac{1}{(q+1)^2} -
\sum_{k=1}^{\lfloor  \frac{q}{2} \rfloor} \tfrac{B_{2k}}{(2k)!}
\tfrac{q! ~(H_q -H_{q-2k+1})}{(q-2k+1)!} -(-1)^q q! \sum_{k\geq \lceil
\frac{q}{2} \rceil +1} \frac{B_{2k}}{(2k)\dots (2k-q-1)} \label{eq:glaisher}
  \end{align}
 \begin{align}
\sum_1 \frac{B_{2k}}{2k} &= \gamma -\half \\
\sum_1 \frac{B_{2k}}{2k-1} &= \zeta'(0) +\gamma +\half = -\half \log
(2\pi) + \gamma +\half\\
\sum_2 \frac{B_{2k}}{2k-2} &= -2\zeta'(-1) +2\zeta'(0) +\gamma
+\tfrac{11}{12} \\ 
\sum_2 \frac{B_{2k}}{2k-3} &= 3\zeta'(-2) -6\zeta'(-1) +3\zeta'(0)
+\gamma +\tfrac{5}{4} \\
\sum_{\lfloor \frac{j}{2} \rfloor +1} \frac{B_{2k}}{2k-j} &= \sum_{i=1}^j
(-1)^{i+j} i {j \choose i} \zeta'(i-j) +\gamma +(H_j -\half
-{\textstyle \sum_{k=1}^{\lfloor \frac{j}{2}
    \rfloor}\frac{B_{2k}}{2k}} ), \qquad (j \geq 0) \\
\sum_{\lfloor \frac{j}{2} \rfloor +1} \frac{B_{2k}}{(2k)(2k-j)} &= \sum_{i=1}^j
(-1)^{i+j} {j-1 \choose i-1} \zeta'(i-j) +\frac{H_j}{j} ,\qquad (j\geq
1) \label{eq:2k(2k-j)} 
  \end{align}
  \begin{align}
\sum B_{2k} \frac{2k-1}{2} &= \zeta(3)-\zeta(2) +\half \label{eq:(2k-1)/2} \\
\sum \frac{B_{2k} (2k-2)}{2k(2k-1)} &= \gamma -\tfrac{3}{2} +\half
\log(2\pi) \label{eq:(2k-2)/(2k(2k-1))}\\
H_n &= \log n +\gamma +\frac{1}{2n} - \sum_1^m \frac{B_{2k}}{2k}
 \frac{1}{n^{2k}} +O\big( \tfrac{1}{n^{2m+1}}\big) \label{eq:H_n} \\
 (i\geq 2:) \qquad H_n^{(i)} &=\zeta(i) -\frac{1}{(i-1) n^{i-1}}
 +\frac{1}{2n^{i}}  - \sum_1^m \frac{B_{2k}}{(2k)!}
 \frac{(i+2k-2)!}{(i-1)!} \frac{1}{n^{i+2k-1}} +O\big(
 \tfrac{1}{n^{i+2m}}\big) \label{eq:H_n^i}
  \end{align}
\end{lem}
\begin{proof}
  The first five lines are the constant terms in the asymptotic
  expansion of $H_n, H_n^{(2)}, H_n^{(3)}, H_n^{(4)}, H_n^{(i)} $.
  The next five lines are the constant terms in
  the asymptotic expansion of $\sum^n k^q \log k $ ($q=0,1,2,3$), the
  first being given by the Stirling formula (for $\log n!$).  These
  are the generalized Glaisher constants \cite{F-03}, see lemma
  \ref{lemma20}.  The third set of lines is obtained recursively by
  partial fraction decomposition from the previous: $\frac{1}{2k-j}=
  \frac{j!}{(2k)...(2k-j)} -\sum_{i=0}^{j-1} \frac{ (-1)^{i+j} \binom
    ji}{2k-i} $; while we used $\frac{1}{(2k)(2k-j)}= -\frac{1}{j}
  (\frac{1}{2k} -\frac{1}{2k-j})$ for (\ref{eq:2k(2k-j)}).  Lines
  (\ref{eq:(2k-1)/2}) and (\ref{eq:(2k-2)/(2k(2k-1))}) are
  miscellaneous linear combinations that we shall use.  The last two
  lines yield the asymptotics for the generalized harmonic numbers
  $\sum_{k=1}^n \frac{1}{k^i}$, with $\zeta(i)$ in (\ref{eq:zeta(i)})
  as it appears in the Euler--Maclaurin formula.
\end{proof}

\begin{lem}
  \begin{align}
\textstyle \sum_1 \frac{B_{2k}}{2k} \frac{(2k)\cdots(2k+i-2)}{(i-1)!}
&=\zeta(i) -\frac{i+1}{2i-2} \\
\textstyle \sum \frac{B_{2k}}{2k} \frac{(2k-1)\cdots(2k+i-3)}{(i-1)!}
&=\zeta(i) -\zeta(i-1) +\frac{1}{(i-1)(i-2)} \\  
\textstyle \sum \frac{B_{2k}}{2k}  \frac{(2k-j)\cdots(2k+i-2-j)}{(i-1)!} &=
\sum_{r=0}^j (-1)^r \binom{j}{r} \zeta(i-r) +\frac{(-1)^{j+1}~j!}{(i-1)\cdots
  (i-j-1)} , \qquad (1 \leq j \leq i-2) \label{eq:zeta-general}\\
&= \sum_{r=0}^{i-1} (-1)^r \binom{j}{r} \zeta(i-r) +c_{j,i}, \qquad
(j\geq i-1, ~i\geq 2) \\
&= \sum_{r=0}^{i-1} \binom{-j-1+r}{r} \zeta(i-r) +c_{j,i}, \qquad
(j\leq -1,~ i\geq 2) 
  \end{align}
with $\zeta(1)$ standing for Euler's $\gamma $ in the last two equations.
As to the constants, they are computed recursively, $ \forall j\in \Z,
i \geq 2$: \\
$ c_{j,i} =c_{j-1,i} -c_{j-1,i-1}$, hence $c_{j,i} =(-1)^{i+1} \Big(
\sum_{r=3}^i (-1)^r c_{j-i-1+r,r} - c_{j-i+2,2} \Big) $, which is
useful for the fourth line; while $ c_{j,i} =c_{j+1,i-1} -c_{j,i-1}$,
hence $ c_{j,i} = \sum_{r=3}^i c_{j+1,r} - c_{j,2} $, is useful for
the fifth line.  The third line can be used for values of $c_{j,i}$
with $ 1\leq j \leq i-2 $.  For any $j$: $c_{j,2} = \frac{j-3}{2}$.  \\
For instance: $ c_{i-2,i} =\frac{(-1)^{i-1}}{i-1}$ $ (
  i\geq 3$, and $c_{0,2}=-\frac{3}{2}$), and $c_{i-1,i} =(-1)^{i-1}
  H_{i-1}$ ($i\geq 2$).\\
\end{lem}
\begin{proof}
These are direct combinatorial consequences of (\ref{eq:zeta(i)}).
\end{proof}
In the lhs of the fourth line ($j\geq i-1$), one can choose to exclude
the few non-zero terms at low values of $k$; in that case $\sum_1$ is
replaced by $\sum_{\lceil \frac{j-i}{2} \rceil}$ and the only change
lies in the constants: $c_{j,2} =\frac{j-3}{2} -\sum_{k=1}^{\lceil
  \frac{j}{2} \rceil -1} B_{2k} \big( 1-\frac{j}{2k} \big)$, and the
recursions remain unchanged.

\section{Stirling numbers and their asymptotics}

For future use, we set 
$$
\fbox{$S_{r,s,t} := {\displaystyle \sum_{s\leq i_1<\dots<i_r \leq t}
  \frac{1}{i_1\cdots i_r}}$}~,
$$
which is 0 for $t < r+s-1$.  Define $S_{0,s,t} :=1$ if $t\geq s-1$
and $S_{r,s,t} :=0$ if $r<0$.  These numbers relate to the Stirling
numbers of the first kind $\big[ {t\atop r} \big]$, defined by
$\sum_{k} \big[ {n \atop k} \big] x^k := x(x+1)\cdots (x+n-1)$, or to
the {\em signed} Stirling numbers $s(t,r)$, defined by $\sum_{k}
s(n,k) x^k := x(x-1)\cdots (x-n+1)$, in the following way: $S_{r,1,t}
= \big[ {t+1 \atop r+1} \big]/t! = (-1)^{r+t} s(t+1,r+1)/t!$. The
generating function for these three versions are:
\begin{align*}
\frac{1}{r!}(\log(1+x))^r  &= \sum_{n\geq 1} s(n,r) \frac{x^n}{n!} \\
\frac{1}{r!}(-\log(1-x))^r  &= \sum_{n\geq 1} \Big[ {n\atop r} \Big]
\frac{x^n}{n!} \\
\frac{1}{r!}(-\log(1-x))^r  &= \sum_{n\geq 1} S_{r-1,1,n-1} \frac{x^n}{n}.
\end{align*}

Here is the relation between Stirling numbers of the first kind and
harmonic numbers.  Denote by $\{r\}=r_1^{i_1}\dots
r_l^{i_l}$ a partition of the integer $r$ into $l$ different parts, ie
$r=i_1 r_1+\dots+i_l r_l$.
Then:
\begin{align} \label{eq:stirling-harmonic}
S_{r,1,n}= (-1)^r \sum_{\{r\}} \prod_{j=1}^l \frac{(-1)^{i_j}}{i_j!}
 \bigg( \frac{H_n^{(r_j)}}{r_j} \bigg)^{i_j}.
\end{align}
The first few cases are $S_{1,1,n}=H_n$ (partition \{1\}=1), and:
\begin{align*}
S_{2,1,n} &= -\half H_n^{(2)} +\half H_n^2 \qquad {\rm partitions}~
\{2\}=2, 1^2.\\
S_{3,1,n} &= \tfrac{1}{3} H_n^{(3)} -\half H_n^{(2)} H_n +\tfrac{1}{6}
H_n^3 \qquad {\rm partitions}~ \{3\}=3, 21, 1^3.\\
S_{4,1,n} &= -\tfrac{1}{4} H_n^{(4)} + \tfrac{1}{3} H_n^{(3)} H_n +\tfrac{1}{8}
(H_n^{(2)})^2 -\tfrac{1}{4} H_n^{(2)} H_n^2 +\tfrac{1}{24} H_n^4 \qquad
{\rm partitions}~ \{4\}=4, 31, 2^2, 21^2, 1^4.
\end{align*}
This can be used, in combination with (\ref{eq:H_n}) and
(\ref{eq:H_n^i}), to compute the asymptotic growth of the Stirling
numbers to arbitrary order.  The terms contributing the most are those
with highest power of $H_n$ ($\sim \log n$): partitions
$\{r\}=1^r,21^{r-2},$ etc. Thus the asymptotic expansion starts as
\begin{align} \label{eq:asymp-stirling}
S_{r,1,n} =\tfrac{1}{r!} (\log n)^r +\tfrac{\gamma}{(r-1)!} (\log
n)^{r-1} +\tfrac{\gamma^2-\zeta(2)}{(r-2)! ~2} (\log n)^{r-2}+\dots  
\qquad {\rm for~ fixed~} r
\end{align}
For an alternative proof of this result and for shedding light on the
decreasing sequence of logarithms, see end of section~\ref{sec:third-sum}. 
Equation (\ref{eq:asymp-stirling}) is the main result of \cite{W-93}
and was strengthened in \cite{H-95}.

As for the asymptotic behaviour when $r$ grows as quickly as $n$, say
for $n-r$ fixed, formula (\ref{eq:stirling-harmonic}) is helpless; but
we can easily find the solution by intuition: $n! S_{n,1,n}=1, ~n!
S_{n-1,1,n}= 1+\dots+n =n(n+1)/2, ~n! S_{n-2,1,n}= \sum_{i_1=1}^n
\sum_{i_2=i_1+1}^n i_1 i_2 =\half $ (square -- diag) $=\half [
n^2(n+1)^2/4 -n(n+1)(2n+1)/6 ] =n(n+1)(3n^2-n-2)/24 \sim n^4/8$.
Similarly: $n! S_{n-3,1,n}= \sum_{i_1=1}^n \sum_{i_2=i_1+1}^n
\sum_{i_3=i_2+1}^n i_1 i_2 i_3 =\frac{1}{3!} $ (cube -- plane) $\sim
\frac{1}{3!} (n(n+1)/2)^3$; and in general we will have $n!
S_{n-k,1,n}\sim \frac{1}{k!} (n(n+1)/2)^k$, that is:
$$
S_{r,1,n} =\frac{1}{n!} \frac{1}{(n-r)!} \Big( \frac{n^2}{2} \Big)^{n-r}
 +\dots \qquad {\rm for}~ n-r ~{\rm const.} 
 $$
 These two asymptotic growths agree with the results of \cite{MW-58}
 obtained by saddle-point evaluation of the generating function
 integral (a method already used by Laplace two centuries ago for
 Stirling numbers of second kind). The same results were re-obtained
 in \cite{KK-91} from recursion equations using the ray method from
 optics.  Formula (\ref{eq:stirling-harmonic}), however, gives as many
 terms as desired for the growth with $r$ fixed.

Formula (\ref{eq:stirling-harmonic}) can be inverted to yield
$$
H_n^{(r)} = (-1)^r r  \sum_{\{r\}} (-1)^{i_1+\dots +i_l}
\frac{(i_1+\dots +i_l-1)!}{i_1!\cdots i_l!} ~S_{r_1,1,n}^{i_1} \cdots
S_{r_l,1,n}^{i_l} .
$$

\section{Asymptotics of sums involving $(\log k)^p/k^q$}

\begin{lem} \label{logk-k}
$ \sum_{k=1}^{n} \frac{\log k}{k} = (\log n) \big[ H_n -\half\log
n -\gamma \big] +\gamma_1 +  \sum_{k=1}^m \frac{B_{2k} H_{2k-1}}{(2k)
  ~n^{2k}} +O(\tfrac 1{n^{2m+1}})$,\\ with $\gamma_1 = - \sum
\frac{B_{2k}}{2k} H_{2k-1}$. 
\end{lem}
\begin{proof} 
  Write the lhs as $\frac{\log n}{n} +\sum_{k=1}^{n-1}
  \frac{\log(k)}{k} $.  For $f(x):=\frac{\log x}{x}$, here are the
  ingredients we need: $\int f(x)= \frac{(\log x)^2}{2}$,
  $f^{(2k-1)}(x)=(2k-1)!\frac{H_{2k-1} -\log x}{x^{2k}}$.  Thus the
  Euler--Maclaurin formula tells us that $\sum_{k=1}^{n-1}
  \frac{\log(k)}{k} = \frac{(\log n)^2}{2} -\half\frac{\log n}{n}+
  \sum \frac{B_{2k}}{(2k) n^{2k}}(H_{2k-1} -\log n) -\sum
  \frac{B_{2k}H_{2k-1}}{(2k)}$.  The $\log n$ terms yield $-(\log n)
  \sum \frac{B_{2k}}{(2k) n^{2k}} $, which estimates $(\log n)[H_n
  -\gamma - \frac{1}{2n}]$ by (\ref{eq:H_n}).  Writing $(\log n) \big[
  H_n -\half\log n -\gamma \big] $ assures us that the remaining terms
  are inverse powers of $n$ (easily tractable under the asymp$_k$
  trick).  The constant $\gamma_1$ is by definition the first
  Stieltjes constant.
\end{proof}

Using the same method of proof, we easily generalise.  
\begin{lem} \label{logkp-k}
$  \sum_{k=1}^{n} \frac{(\log k)^2}{k} = \frac{(\log n)^3}{3} 
+\gamma_2 +\half\frac{(\log n)^2}{n} -\sum_{k=1}^m \frac{B_{2k}}{2k}
\frac{(\log n)^2 -2 H_{2k-1}(\log n) +H^2_{2k-1}-H^{(2)}_{2k-1}}{n^{2k}}
+O(\tfrac {(\log n)^{2}}{n^{2m+1}})$,\qquad with $\gamma_2 = \sum \frac{B_{2k}
  (H^2_{2k-1} -H^{(2)}_{2k-1})}{2k}$.  Similarly, for $p\geq 0$: 
$$ 
\sum_{k=1}^{n} \frac{(\log k)^p}{k} = \tfrac{1}{p+1} (\log n)^{p+1} 
+\gamma_p +\half\frac{(\log n)^p}{n} -\sum_{k=1}^m  
\frac{\sum_{r=0}^p d_{p,k,r}(\log n)^{p-r} }{n^{2k}}
+O(\tfrac {(\log n)^{p}}{n^{2m+1}}),
$$
with $ d_{p,k,r}= \frac{B_{2k}}{2k} \frac{(-1)^r p!}{(p-r)!} ~S_{r,1,2k-1}$
~~and $\gamma_p = (-1)^p ~p! {\displaystyle \sum_{k>p/2}} 
\frac{B_{2k}}{2k} ~S_{p,1,2k-1}$.
\end{lem}
Here, $\gamma_p$ is the $p$-th Stieltjes constant by definition.  For
$p=0$ we have $S_{0,s,t}=1$ and $\gamma_0=\gamma$ (though one needs to add
$\half$ to $\gamma_p$ in this case).  Note that in the last sum, we
could drop the requirement $k> p/2$, as $S_{p,1,2k-1}$ vanishes for
$k=1,\dots,p/2$. We will neglect this in future.

Note that the expressions we find for $\gamma_p$ are exactly the same
as those one finds when directly analytically expanding the zeta
function via the Euler--Maclaurin formula. One uses $f(x):=x^{-s}$ with
$f^{(2k-1)}(x) =-s(s+1)\cdots (s+2k-2) ~x^{-s-2k+1} = -
(2k-1)! ~\sum_{p=0}^{2k-1} S_{p,1,2k-1} (s-1)^p  ~x^{-s-2k+1} $:
\begin{align*}
\zeta(s)= \lim_{N\to\infty} \sum_{n=1}^N \frac{1}{n^s}
& = \frac{1}{s-1} +\half +\sum_{k=1}^M \frac{B_{2k}}{(2k)!}~
s(s+1)\cdots (s+2k-2) + {\rm error}(s,M) \\ &=  \frac{1}{s-1} +
\sum_{p\geq 0} \Big( \underbrace{ \sum_{k>p/2} \frac{B_{2k}}{2k}
  S_{p,1,2k-1} }_{\frac{(-1)^p}{p!} ~\mbox{\large $\gamma_p$}} \Big) (s-1)^p .
\end{align*}

We also easily generalise in another direction:
\begin{lem} \label{logk-k2}
$  \sum_{k=1}^{n} \frac{\log k}{k^2} = -\zeta'(2) + (\log n)
\big[ H^{(2)}_n -\zeta(2) \big] -\frac{1}{n} +\sum_{k=1}^m
\frac{B_{2k} (H_{2k}-1)}{n^{2k+1}} +O(\tfrac 1{n^{2m+2}})$, \\
with $-\zeta'(2) = 1-\sum B_{2k} (H_{2k}-1)$.  Similarly, for $q\geq 2$: 
$$ 
 \sum_{k=1}^{n} \frac{\log k}{k^q} = -\zeta'(q) + (\log n)
\big[ H^{(q)}_n -\zeta(q) \big] -\tfrac{1}{(q-1)^2 n^{q-1}}
 +\sum_{k=1}^m \tfrac{B_{2k}}{(2k)!} \tfrac{(q+2k-2)!}{(q-1)!}
\tfrac{H_{q+2k-2}-H_{q-1}}{n^{q+2k-1}} +O(\tfrac 1{n^{q+2m}}),
$$
with $-\zeta'(q) =\frac{1}{(q-1)^2} -\sum
\frac{B_{2k}}{(2k)!} \frac{(2k+q-2)!}{(q-1)!} (H_{2k+q-2}-H_{q-1})$.
\end{lem}

The meta-generalisation regroups the two previous results:
\begin{lem} \label{logkp-kq}  For $q\geq 2$ and $p\geq 0$:
$$
\sum_{k=1}^{n} \frac{(\log k)^p}{k^q} = (-1)^p \zeta^{(p)}(q)
-\frac{\sum_{r=0}^p \frac{p!/(p-r)!}{(q-1)^{r+1}} (\log
  n)^{p-r}}{n^{q-1}} +\half \frac{(\log n)^p}{n^q} - \sum_{k=1}^m
\frac{\sum_{r=0}^p d_{p,q,k,r} (\log n)^{p-r}}{n^{q+2k-1}} +O(\tfrac
{(\log n)^{p}}{n^{q+2m}})
$$
with $d_{p,q,k,r} = \frac{B_{2k}}{(2k)!} \frac{(2k+q-2)!}{(q-1)!}
\frac{(-1)^r p!}{(p-r)!}  S_{r,q,2k+q-2} $, \\
and $(-1)^p \zeta^{(p)}(q) = \frac{p!}{(q-1)^{p+1}} +(-1)^p p!  \sum_1
\frac{B_{2k}}{(2k)!}  \frac{(2k+q-2)!}{(q-1)!}  S_{p,q,2k+q-2} $.
\end{lem}
\begin{proof}
  Write the lhs as $\frac{(\log n)^p}{n^q} +\sum_{k=1}^{n-1}
  \frac{(\log k)^p}{k^q} $.  For $f(x):=\frac{(\log x)^p}{x^q}$, here
  are the ingredients we need: $\int f(x)= -\frac{1}{(q-1)^{r+1}
    x^{q-1}} \sum_{r=0}^p \frac{p! (\log x)^{p-r}}{(p-r)!}$,
  $f^{(i)}(x) = \frac{1}{x^{q+i}} \frac{(q+i-1)!}{(q-1)!} \sum_{r=0}^p
  \frac{(-1)^{r+i} p!}{(p-r)!} S_{r,q,q+i-1} (\log x)^{p-r}$ (NB:
  $S_{...} =0$ for $r>i$).  Now simply apply the Euler--Maclaurin
  formula.
\end{proof}

\paragraph{Application to numerics of the $\zeta$ function.}
Note that this formula, together with the asymp$_k$ trick, allows a
very rapid numerical computation of $\zeta^{(p)}(q)$ (for positive
integer $q$), much more efficient than current mathematical softwares.
Lemmas \ref{kq-logkp} or \ref{logkp-(n-k)} provide a formula for
negative $q$.

\section{Asymptotics of sums involving $k^q(\log k)^p$}

\begin{lem}  For $p\geq 1$ we have:
$$
\sum_{k=1}^n (\log k)^p = n \sum_{r=0}^p \tfrac{(-1)^r p!}{(p-r)!}
(\log n)^{p-r} +\half (\log n)^p + {\rm const} - \sum_{k=1}^m
\frac{\tfrac{B_{2k}}{2k(2k-1)}  \sum_{r=1}^p \tfrac{(-1)^r
    p!}{(p-r)!} S_{r-1,1,2k-2} (\log n)^{p-r} }{n^{2k-1}}
+\CO(\tfrac {(\log n)^{p-1}}{n^{2m}}), 
$$
with {\rm const} $=(-1)^p \zeta^{(p)}(0) =(-1)^{p} p! \big( -1+\sum_1
\frac{B_{2k}}{2k(2k-1)}  S_{p-1,1,2k-2} \big) $.
\end{lem}

\begin{lem} \label{lemma20}  For $q\geq 0$ we have:
  \begin{align*}
\sum_{k=1}^{n} k^q \log k &=\frac{n^{q+1}}{q+1}\big[ (\log n)
-\tfrac{1}{q+1} \big] +\half n^q \log n
+\sum_{k=1}^{\lceil \frac{q}{2} \rceil}
\tfrac{B_{2k}}{(2k)!}\tfrac{q! }{(q-2k+1)!}  \big[ (\log n) 
+(H_q -H_{q-2k+1)} \big] ~n^{q-2k+1}\\
& + {\rm const}  +(-1)^q q! \sum_{k= \lceil
  \frac{q}{2} \rceil +1}^m \tfrac{B_{2k}}{(2k)\dots (2k-q-1)}
\tfrac{1}{n^{2k-q-1}} +O(\tfrac 1{n^{2m-q}})
  \end{align*}
  with ${\rm const} = -\zeta'(-q) +\frac{B_{q+1}}{q+1} H_q
  =\frac{1}{(q+1)^2} -\sum_{k=1}^{\lfloor \frac{q}{2} \rfloor}
  \tfrac{B_{2k}}{(2k)!}  \tfrac{q! (H_q -H_{q-2k+1})}{(q-2k+1)!}
  -(-1)^q q! \sum_{k\geq \lceil \frac{q}{2} \rceil +1}
  \tfrac{B_{2k}}{(2k)\dots (2k-q-1)} $ being the generalized Glaisher
  constant of (\ref{eq:glaisher}).  For odd $q$, it is understood that
  the last term of the sum $\sum_{k=1}^{\lceil \frac{q}{2} \rceil}
  \dots$ (with $k=\lceil \frac{q}{2} \rceil$ and independent of $n$),
  which equals $\frac{B_{q+1}}{q+1} H_q$, should not be counted as it
  is already counted in {\rm const}.
\end{lem}

Again, the meta-generalisation regroups the two previous results:
\begin{lem} \label{kq-logkp} For $q\geq 0$ and $p\geq 1$ we have:
  \begin{align*}
\sum_{k=1}^{n} k^q (\log k)^p &= n^{q+1} \sum_{r=0}^p \tfrac{(-1)^r
p!}{(p-r)!} \tfrac{(\log n)^{p-r}}{(q+1)^{r+1}} +\half n^q 
(\log n)^p + \sum_{k=1}^{\lceil \frac{q}{2} \rceil} \Big( \sum_{r=1}^p
c_{p,q,k,r} (\log n)^{p-r} \Big) ~n^{q-2k+1} +  \\ 
&~~ + {\rm const} + \sum_{k= \lceil \frac{q}{2} \rceil +1}^m \frac{\sum_{r=1}^p
  d_{p,q,k,r} (\log n)^{p-r}}{n^{2k-q-1}} +O(\tfrac {(\log n)^{p-1}}{n^{2m-q}})
  \end{align*}
  with $c_{p,q,k,r} := \tfrac{B_{2k}}{(2k)!}\tfrac{q!}{(q-2k+1)!}
  \tfrac{p!}{(p-r)!} S_{r,q-2k+2,q}$  \\
  and $ d_{p,q,k,r} := \tfrac{B_{2k} (-1)^{r+q+1}}{(2k)\cdots(2k-q-1)}
  \tfrac{p!q!}{(p-r)!} \sum_{j=0}^q (-1)^j S_{j,1,q}
  S_{r-j-1,1,2k-q-2}$  \\
  and ${\rm const} = (-1)^p \zeta^{(p)} (-q) +\frac{B_{q+1}}{q+1} p!
  S_{p,1,q} =\frac{(-1)^{p+1} p!}{(q+1)^{p+1}} -\sum_{k=1}^{\lfloor
    \frac{q}{2} \rfloor}
  c_{p,q,k,p} - \sum_{k\geq \lceil \frac{q}{2} \rceil +1} d_{p,q,k,p}$.\\
  For odd $q$, it is understood that the last term of the sum
  $\sum_{k=1}^{\lceil \frac{q}{2} \rceil} \dots$ (with $k=\lceil
  \frac{q}{2} \rceil$ and $r=p$), which is constant and equals
  $\frac{B_{q+1}}{q+1} p! S_{p,1,q}$, should not be counted as it is
  already counted in {\rm const}.
\end{lem}

\begin{proof}
  Write the lhs as $ n^q (\log n)^p +\sum_{k=1}^{n-1}
   k^q (\log k)^p $.  For $f(x):= x^q (\log x)^p $, here
  are the ingredients we need: $\int f(x)= \frac{x^{q+1}}{(q+1)^{r+1}
    } \sum_{r=0}^p \frac{ (-1)^p p!}{(p-r)!} (\log x)^{p-r}$ and
 $$
f^{(i)}(x) = \left\{     \begin{array}{l}
   \frac{x^{q-i}}{(q-i)!} \sum_{r=0}^p \frac{p!}{(p-r)!} (\log
   x)^{p-r} q! S_{r,q-i+1,q} \qquad {\rm for} ~i \leq q, ~~({\rm NB}:
   S_{...} =0 ~{\rm for}~  r>i)\\ 
 \frac{(i-q-1)!}{x^{i-q}} \sum_{r=0}^p
  \frac{(-1)^{r+q+i} p!}{(p-r)!} (\log x)^{p-r} q! \sum_{j=0}^q (-1)^j
  S_{j,1,q} S_{r-j-1,1,i-q-1} \qquad {\rm for}~ i>q,
    \end{array}  \right.
$$
Now simply apply the Euler--Maclaurin formula.
\end{proof}

\section{Asymptotics of sums involving $(\log k)^p/(n-k)^q$}
\label{sec:third-sum}
\begin{lem} \label{logk-(n-k)}
$  \sum_{k=1}^{n-1} \frac{\log k}{n-k} = (\log n)^2 +\gamma (\log n)
-\zeta(2) +\sum_{k=1}^m \big( \frac{(-1)^k B_k}{k^2} -\zeta'(1-k) \big)
\frac{1}{n^k} +O(\tfrac 1{n^{m+1}})$.
\end{lem}
\begin{proof} 
  Write the lhs as $\log (n-1) +\sum_{k=1}^{n-2} \frac{\log(k)}{n-k} $.
  The integral of $f(x):=\frac{\log x}{n-x}$ is: $-\log (x)
  \log(1-\frac{x}{n}) -\li(\frac{x}{n})$, thus\footnote{We have used:
    $\li(1-x) +\li(x) = -\log(x) \log(1-x) +\zeta(2)$ (proof by
    derivation).}  $\int_1^{n-1} f(x) = (\log n)^2 -\zeta(2) +2
  \li(\frac{1}{n})$.  Note also that $f^{(i)}(x) = i! \frac{\log
    x}{(n-x)^{i+1}} +i! \sum_{r=1}^i \frac{(-1)^{r+1} /r}{
    (n-x)^{i+1-r} x^r}$, so that $ [f^{(2k-1)}(x)]^{n-1}_1 = (2k-1)!
  \big[ \log (n-1) + \sum_{r=1}^{2k-1} \frac{(-1)^{r+1}}{r} \big(
  \frac{1}{(n-1)^r} - \frac{1}{(n-1)^{2k-r}} \big) \big] $.  In total:
{\scriptsize  
  \begin{align*}
 \sum_{k=1}^{n-1} \frac{\log(k)}{n-k} &= (\log n)^2 -\zeta(2)
  +2~\li\big(\frac{1}{n}\big) + \gamma \log(n-1) +\frac{1}{(n-1)} \sum_1
  \frac{B_{2k}}{2k} \Big( 1 -\frac{1}{2k-1} \Big) - \frac{1}{(n-1)^2}
  \sum_2 \frac{B_{2k}}{2k} \Big( \frac{1}{2} -\frac{1}{2k-2} \Big) +
  \dots\\
&~~ + \frac{ (-1)^{j+1}}{(n-1)^j} \sum_{\lfloor \frac{j}{2}
    \rfloor +1} \frac{B_{2k}}{2k} \Big( \frac{1}{j} -\frac{1}{2k-j}
  \Big) +\dots  
  \end{align*}
} Now use (\ref{eq:2k(2k-j)}) as well as the expansions $\log
(n-1)=\log n -(\frac{1}{n} +\frac{1}{n^2}+\dots)$ and
$\frac{1}{(n-1)^j} = \sum_{i\geq 0} \binom{j+i-1}{i-1}
\frac{1}{n^{j+i}}$.  There are nice cancellations so that only
$-\zeta'(1-k)$ survives at power $\frac{1}{n^k}$.
\end{proof}

\begin{lem} \label{logkp-(n-k)} For $p \geq 1$ we have:
  \begin{align*}
 \sum_{k=1}^{n-1} \frac{(\log k)^p}{n-k} &= (\log n)^{p+1} +\gamma
 (\log n)^{p} + \sum_{r=1}^{p} c_{p,r} (\log n)^{p-r}   \\
&~~ -\sum_{k=1}^m
 \Big( \tfrac{(-1)^k B_k}{k^2} \sum_{r=1}^{p} d_{p,r,k} (\log n)^{p-r} 
 -(-1)^p \zeta^{(p)}(1-k) \Big) \frac{1}{n^k} +\CO(\tfrac {(\log
 n)^{p-1}}{n^{m+1}})
  \end{align*}
with $c_{p,r} := \tfrac{(-1)^r p!}{(p-r)!} \zeta(r+1)$ and 
$ d_{p,r,k} := \tfrac{(-1)^r p!}{(p-r)!} S_{r-1,1,k-1}$.
Hence, the constant term is $ (-1)^p p! \zeta(p+1) $.
\end{lem}
\begin{proof} 
  Write the lhs as $\log^p(n-1) +\sum_{k=1}^{n-2}
  \frac{(\log k)^p}{n-k} $.  With $f(x):=\frac{(\log
    x)^p}{n-x}$, we have:
{\scriptsize
  \begin{align*} 
\int f(x) &= -\log(1-\tfrac{x}{n}) (\log x)^p + \sum_{r=1}^{p}
\tfrac{(-1)^r p!}{(p-r)!} \Li_{r+1}(\tfrac xn) (\log x)^{p-r} \\
\int_1^{n-1} f(x) &= \log(n) \log^p(n-1) + \sum_{r=1}^{p}
\tfrac{(-1)^r p!}{(p-r)!} \Li_{r+1}(1-\tfrac 1n) \log^{p-r}(n-1)
-(-1)^p p! \Li_{p+1}(\tfrac 1n)  \\
&= \log(n) \log^p(n-1) - \sum_{r=1}^{p}
\tfrac{(-1)^r p!}{(p-r)!} \big( \mathfrak{S}_{1,r}(\tfrac 1n)
-\zeta(r+1) \big) (\log n)^{p-r} -(-1)^p p! \Li_{p+1}(\tfrac 1n)  \\
 f^{(i)}(x) &= i! \frac{(\log x)^p}{(n-x)^{i+1}} +i! \sum_{j=1}^i
 \frac{(-1)^j/j}{(n-x)^{i+1-j} x^j} \sum_{r=1}^{p} \tfrac{(-1)^r
 p!}{(p-r)!} S_{r-1,1,j-1} (\log x)^{p-r} \\
[f^{(2k-1)}(x)]^{n-1}_1 &= (2k-1)! \log^p(n-1) +(2k-1)!
 \sum_{j=1}^{2k-1} \frac{(-1)^j}{(n-1)^j} \Big( \sum_{r=1}^{p} \tfrac{(-1)^r
 p!}{(p-r)!}  \log^{p-r}(n-1) \frac{S_{r-1,1,j-1}}{j} -(-1)^p
 \frac{S_{p-1,1,2k-j-1}}{2k-j} \Big)
  \end{align*}
}
where we used the notation $ \mathfrak{S}_{1,r}(x) := \sum_{n\geq 1}
\frac{S_{r-1,1,n-1}}{n^2} x^n $ from section~\ref{sec:polylogs}.  In total:
{\scriptsize
  \begin{align*}
    \sum_{k=1}^{n-1}
  \frac{(\log k)^p}{n-k}  &= \int_1^{n-1} f(x) + \half \log^p(n-1) +
(\gamma -\half) \log^p(n-1)\\
&~~~  - \frac{1}{(n-1)} \sum_{k\geq 1} \frac{B_{2k}}{2k}
\Big(  \sum_{r=1}^{p} \tfrac{(-1)^r p!}{(p-r)!} \log^{p-r}(n-1) 
  \frac{S_{r-1,1,0}}{1} -(-1)^p \frac{S_{p-1,1,2k-2}}{2k-1} \Big) \\
&~~~ \pm \dots + \frac{(-1)^j}{(n-1)^j} \sum_{k\geq \lfloor \frac{j}{2}
  \rfloor +1} \frac{B_{2k}}{2k} \Big( \sum_{r=1}^{p} \tfrac{(-1)^r p!}{(p-r)!}
  \log^{p-r}(n-1) \frac{S_{r-1,1,j-1}}{j} -(-1)^p
  \frac{S_{p-1,1,2k-j-1}}{2k-j} \Big) +\dots
  \end{align*}
}
The remainder of the proof are nice cancellations, which are
impossible to prove in the general case; we exhibit here the case
$p=2$ as a pattern for all other cases.  For \underline{$p=2$} we have:
{\scriptsize
  \begin{align*}
    \sum_{k=1}^{n-1} \frac{(\log k)^2}{n-k}  &= \int_1^{n-1} f(x) +
    \half \log^2(n-1) + \sum_1 \frac{B_{2k}}{2k} [f^{(2k-1)}(x)]^{n-1}_1\\
&= (\log n)^3 +2(\log n) [\mathfrak{S}_{1,1}-\zeta(2)]
    -2[\mathfrak{S}_{1,2}(\tfrac 1n) -\zeta(3) 
    +\Li_3(\tfrac 1n)] + \half \log^2(n-1) +\\
&~~~ + (\gamma -\half) \log^2(n-1) +\frac{2}{(n-1)} \sum_{1} \frac{B_{2k}}{2k} 
\Big( \log(n-1) +\frac{H_{2k-2}}{2k-1} \Big) - \frac{2}{(n-1)^2}
    \sum_{2} \frac{B_{2k}}{2k} \Big( \frac {\log(n-1)}{2}
    -\frac{H_1}{2}+\frac{H_{2k-3}}{2k-2} \Big) \\
&~~~ \pm \dots +\frac{2(-1)^j}{(n-1)^j} \sum_{\lfloor \frac{j}{2}
  \rfloor +1} \frac{B_{2k}}{2k} \Big( \frac {\log(n-1)}{j}
    -\frac{H_{j-1}}{j} +\frac{H_{2k-j-1}}{2k-j} \Big) +\dots
  \end{align*}
} 
Now replace $\sum_{\lfloor \frac{j}{2} \rfloor +1} \frac{B_{2k}}{2k}
$ by $\gamma-\half -\sum_{r=1}^{\lfloor \frac{j}{2} \rfloor}
\frac{B_{2r}}{2r} $ and set $ h_j := 2 \sum_{\lfloor \frac{j}{2}
  \rfloor +1} \frac{B_{2k}}{2k} \frac{H_{2k-j-1}}{2k-j}$.  The
($\gamma-\half$) will cancel out due to $\log (1-\frac{1}{n}) = -\log
(1-\frac{1}{1-n}) =\sum_{i\geq 1} \frac{(-1)^i}{i (n-1)^i}$ (and the
squared version of it).  Now use (\ref{eq:S11}) and (\ref{eq:S12}) to
simplify the rest and arrive at
{\scriptsize
  \begin{align*}
    \sum_{k=1}^{n-1} \frac{(\log k)^2}{n-k}  &= (\log n)^3
    +\gamma(\log n)^2 -2\zeta(2) (\log n) +2\zeta(3) + \sum_{k\geq 1}
    \Big[ \frac{(-1)^k B_k}{k^2} 2(\log n -H_{k-1}) -\frac{2}{k^3}
    -\sum_{r=1}^k (-1)^r \binom{k-1}{r-1} h_r \Big] \frac{1}{n^k}.
  \end{align*}
} 
Use {\small $ \zeta''(-q) =-\frac{2}{(q+1)^3} -\frac{B_{q+1}}{q+1}
2S_{2,1,q} -\sum_{k=1}^{\lfloor \frac{q}{2} \rfloor} \frac{B_{2k}}{(2k)!}
\frac{q!}{(q-2k+1)!} 2S_{2,q-2k+2,q} - \sum_{k\geq \lceil \frac{q}{2} \rceil
  +1} \frac{B_{2k} (-1)^q 2q!}{(2k)\cdots (2k-q-1)} (H_{2k-q-2} -H_q)
$} from lemma \ref{kq-logkp} as well as the partial fraction
decomposition $ \frac{1}{2k(2k-j)} =\sum_{r=1}^j \frac{(j-1)!}{(j-r)!}
\frac{1}{(2k)\cdots (2k-r)}$ to show that 
{\small
$$  
 h_j :=  ~~2 \sum_{k\geq \lfloor \frac{j}{2} \rfloor +1} \frac{B_{2k}}{2k}
\frac{H_{2k-j-1}}{2k-j} \quad = -\sum_{r=1}^j (-1)^{r} \binom{j-1}{r-1}
\big[\zeta''(1-r)+\tfrac{2}{r^3} \big] 
$$
}
or equivalently: $ -\sum_{r=1}^k (-1)^r \binom{k-1}{r-1} h_r
=\zeta''(1-k) +\frac{2}{k^3} $.  This completes the proof for $p=2$.
The proofs for $p\geq 3$ run similarly.
\end{proof}

\begin{lem} \label{logk-(n-k)q} For $q \geq 2$ we have:
$$
 \sum_{k=1}^{n-1} \frac{\log k}{(n-k)^q} = \zeta(q) \log n
- \sum_{i=1}^{q-2} \tfrac{\zeta(q-i)}{i~n^i}  -\tfrac{2 \log n
  +C_{q-1}}{(q-1)~ n^{q-1}} +\sum_{k=1}^m \textstyle \Big( \frac{(-1)^k
  B_k}{k(k+q-1)} -\binom{k+q-2}{q-1} \zeta'(1-k) \Big)
\frac{1}{n^{k+q-1}} +O(\tfrac 1{n^{m+q}})
$$
with $C_{q} := \gamma- H_q + \frac{2}{q}$.  \qquad Hence, there is no
constant term. 
\end{lem}
\begin{proof} 
  Write the lhs as $\log (n-1) +\sum_{k=1}^{n-2}
  \frac{\log(k)}{(n-k)^q} $.  With $f(x):=\frac{\log
    x}{(n-x)^q}$, we have:
{\scriptsize
  \begin{align*} 
\int f(x) &= -\frac{1}{q-1} (\log x) \Big( \frac{1}{x^{q-1}}
  -\frac{1}{(x-n)^{q-1}} \Big) + \frac{1}{q-1} \frac{\log
    (n-x)}{n^{q-1}} - \sum_{i=1}^{q-2} \frac{1}{i(q-1) ~n^{q-i-1}
    (n-x)^i} \\  
\int_1^{n-1} f(x) &= \frac{1}{q-1} \log (n-1) +
  \sum_{i=1}^{q-2} \frac{1}{(q-1) i ~n^{q-1-i}} \Big(
  \frac{1}{(n-1)^i} -1 \Big) - \frac{2}{q-1} \frac{\log
    (n-1)}{n^{q-1}} \\
 f^{(i)}(x) &=
  \frac{(q-1+i)!}{(q-1)!}\frac{\log x}{(n-x)^{q+i}} +i! \sum_{r=1}^i
  \frac{ \frac{(-1)^{r+1}}{r} \binom{i-r+q-1}{q-1} }{ (n-x)^{q+i-r}
    x^r} \\
[f^{(2k-1)}(x)]^{n-1}_1 &=
  \frac{(q+2k-2)!}{(q-1)!} \log(n-1)+ (2k-1)!  \sum_{r=1}^{2k-1}
  \frac{(-1)^{r+1}}{r} \binom{2k-r+q-2}{q-1} \Big( \frac{1}{(n-1)^r} -
  \frac{1}{(n-1)^{2k-r+q-1}} \Big) .  
  \end{align*}
}
In total: 
{\scriptsize
  \begin{align*}
    \sum_{k=1}^{n-1}
  \frac{\log k}{(n-k)^q}  &= \int_1^{n-1} f(x) + \half \log (n-1) + \log
  (n-1) \sum_1 \frac{B_{2k}}{(2k)!}  \frac{(q+2k-2)!}{(q-1)!} \\
&~~  + \frac{1}{n-1} \sum_1 \frac{B_{2k}}{2k} \binom{2k+q-3}{q-1} -
  \frac{1}{(n-1)^2} \frac{1}{2} \sum_1 \frac{B_{2k}}{2k} \binom{2k+q-4}{q-1}
  +\dots + \frac{1}{(n-1)^{q-1}} \frac{(-1)^{q}}{q-1} \sum_1
  \frac{B_{2k}}{2k} \binom{2k-1}{q-1} \\ 
&~~  + \frac{1}{(n-1)^q} \Bigg[ \frac{(-1)^{q+1}}{q} \sum_1
  \frac{B_{2k}}{2k} \binom{2k-2}{q-1} - \sum_1
  \frac{B_{2k}}{(2k)(2k-1)} \Bigg] +\dots \\ 
&~~  + \frac{(-1)^{j+1}}{(n-1)^{q+j}} \Bigg[ \frac{(-1)^q}{q+j}
  \sum_{\lfloor \frac{j+1}{2} \rfloor +1} \frac{B_{2k}}{2k}
  \binom{2k-j-2}{q-1} - \binom{j+q-1}{q-1} \sum_{\lfloor \frac{j+1}{2}
    \rfloor +1} \frac{B_{2k}}{(2k)(2k-j-1)} \Bigg] +\dots 
  \end{align*}
}
Now use (\ref{eq:2k(2k-j)}), (\ref{eq:zeta-general}), as well as the
expansions $\log (n-1)=\log n -(\frac{1}{n} +\frac{1}{n^2}+\dots)$ and
$\frac{1}{(n-1)^j} = \sum_{i\geq 0} \binom{j+i-1}{i-1}
\frac{1}{n^{j+i}}$.  There are nice cancellations so that only
$-\zeta'(1-k)$ survives at power $\frac{1}{n^k}$.
\end{proof}

Note that without going through the proof, one can empirically
determine the values of the $C_q$ just using the asymp$_k$ trick: one
first uses the trick to quickly determine the 30 first values of
$C_q$, then uses it again to determine the asymptotic growth of those
values up to $\CO(\frac{1}{n^6})$ and recognizes the growth of
harmonic numbers.

Again, the meta-generalisation regroups the two previous results:
\begin{lem} \label{logkp-(n-k)q}
   For $q \geq 2$ and $p\geq 1$ we have:
   \begin{align*}
 \sum_{k=1}^{n-1} \frac{(\log k)^p}{(n-k)^q} &= \zeta(q) (\log n)^p
+ \sum_{i=1}^{q-2} \tfrac{\zeta(q-i)}{i~n^i} \sum_{r=1}^i c_{p,i,r}
(\log n)^{p-r}  +\tfrac{1}{(q-1)~ n^{q-1}} \sum_{r=0}^p
d_{p,q,r} (\log n)^{p-r} \\
&~~ -\sum_{k=1}^m \bigg( \tfrac{(-1)^k
  B_k}{k(k+q-1)} \sum_{r=1}^p c_{p,k+q-1,r}
(\log n)^{p-r} -(-1)^p \tbinom{k+q-2}{q-1} \zeta^{(p)}(1-k) \bigg)
\frac{1}{n^{k+q-1}} +O(\tfrac {(\log n)^{p-1}}{n^{m+q}})     
   \end{align*}
with~ $c_{p,i,r} := \frac{(-r)^r p!}{(p-r)!} S_{r-1,1,i-1} $\\
and ~$d_{p,q,r} := \frac{(-r)^r p!}{(p-r)!} \times 
\left\{    \begin{array}{l}
D_{r,q} -(p-r) S_{r,1,q-2} +S_{r-1,1,q-2}~ \gamma +\sum_{s=2}^r
S_{r-s,1,q-2} ~\zeta(s) ~~{\rm for}~ r=0,\dots,q-1\\
 \sum_{s=r-q+2}^r S_{r-s,1,q-2}
 \zeta(s) ~~~{\rm for}~~ r=q,\dots,p ~~({\rm in ~case}~ p\geq q),
    \end{array} \right. $
wherein $D_{r,q}$ are  the rational numbers \scriptsize{
$$
D_{r,q} := \sum_{r=0}^{p-2} \sum_{j=1}^{q-p+r}
\frac{S_{r,1,j-1}}{j}  \bigg( \frac{S_{p-2-r,1,q-j-2}}{j}
-\sum_{r=1}^{q-j-2} \frac{\sum_{i=0}^{r-1} (-1)^i S_{i,1,r-1}
  S_{p-3-r-i,1,q-j-2-r}}{r(r+j)}  \bigg) +  \frac{S_{p-1,1,q-2}}{q-1}
-S_{p,1,q-2}.
$$
}  \normalsize
In particular, there is no constant term in the asymptotic expansion.
\end{lem}

\begin{proof}
  Write the lhs as $\log^p(n-1) +\sum_{k=1}^{n-2}
  \frac{(\log k)^p}{(n-k)^q} $.  With $f(x):=\frac{(\log
    x)^p}{(n-x)^q}$, we have:
{\scriptsize
  \begin{align*} 
    \int f(x) &= \frac{(\log x)^p}{(q-1)!} \sum_{j=0}^{q-2}
    \tbinom{q-1}{j} \frac{x^{q-j-1}}{n^{q-1} (n-x)^{q-1-j}}\\
  &~~+\sum_{k=1}^p \frac{(-1)^k p!}{(p-k)!} (\log x)^{p-k} \bigg(
    \sum_{j=1}^{q-2} \Big( \sum_{r=0}^{j-1} \frac{S_{k-1,r} r!
      (q-r-3)!} {(j-r-1)! (q-j+r)!} \Big) \frac{x^{q-1-j}} {n^{q-1}
    (n-x)^{q-1-j}} +\sum_{j=1}^k \frac{S_{j-1,q-2}}{(q-1) n^{q-1}}
    \Li_{k-j+1}(\tfrac xn) \bigg)\\
\int_1^{n-1} f(x) &= \bigg[\frac{p!}{(q-1) n^{q-1}} \bigg(
    \sum_{l=0}^p \frac{(-1)^l}{(p-l)!}\Big( -S_{l,1,q-2}
    +\sum_{j=1}^{q-l-1} d_{q,l,j} n^j +S_{l-1,1,q-2} (\log n) \Big) \\
   &~~~~ +\sum_{l=0}^{p-2} \Big( \sum_{r=0}^l \frac{(-1)^{r+p}}{r!}
    S_{l-r,1,q-2} \log^r(n-1) \Big) \Li_{p-l}(1-\tfrac 1n)\bigg) \bigg]
  -\bigg[(-1)^p p! \sum_{j=1}^p \frac{S_{j-1,1,q-2}}{(q-1) n^{q-1}}
    \Li_{p-j+1}(\tfrac 1n) \bigg]
 \end{align*} 
 \begin{align*} 
 f^{(i)}(x) &= \frac{(q-1+i)!}{(q-1)!} \frac{(\log x)^p}{(n-x)^{q+i}}
   - \sum_{r=1}^p \frac{p!(\log x)^{p-r}}{(p-r)! (q-1)!} \sum_{j=1}^{i+1-r}
     \frac{(-1)^j i!(q+i-j+r)!}{(j+r-1) (i-j-r+1)!}
    \frac{ S_{r-1,1,j+r-2}}{(n-x)^{q+i-r-j+1} x^{r+j-1}} \\ 
 [f^{(2k-1)}(x)]^{n-1}_1 &= \frac{(q+2k-2)! \log^p(n-1)}{(q-1)!}
    -\frac{(2k-1)!}{(q-1)!} \sum_{j=1}^{2k-1} \frac{1}{(n-1)^j} \sum_{r=1}^j 
   \frac{(-1)^{j-r} p!\log^{p-r}(n-1)}{(p-r)!} 
    \frac{ (q+2k-j-2)!}{ j ~(2k-1-j)!} S_{r-1,1,j-1} \\
   & \hspace{3.5cm} +\frac{(2k-1)!}{(q-1)!}
    \sum_{j=0}^{2k-p-1}  \frac{1}{(n-1)^{q+j}} \frac{(-1)^{p+j} p!
    (q+j-1)!}{(2k-j-1) ~j!} S_{p-1,1,2k-j-2}
  \end{align*}
} The total expression for $ \sum_{k=1}^{n-1} \frac{(\log k)^p}{(n-k)^q}$
from the Euler--Maclaurin formula is too messy to write out.  As
usual, cancellations will be hard at work and the result will boil
down to the rhs in the lemma.  The closed expression for the $D_{r,q}$
was particularly hard to find (empirically).
\end{proof}

\subsubsection*{Application to the asymptotics of Stirling numbers.}
Had we not known the asymptotic growth of Stirling numbers
(\ref{eq:asymp-stirling}), we could easily find it by induction from
the leading terms in lemma \ref{logkp-k} and \ref{logkp-(n-k)}.
Assuming the empirical result (via the asymp$_k$ trick) that the
coefficient of $x^n$ in $(-\log(1-x))^p$ has leading behaviour
$\sim p(\log n)^{p-1}/n$, we prove:
\begin{align*}
  (p+1)! \sum_{n\geq 1} S_{p,1,n-1} \frac{x^n}{n} &=
  (-\log(1-x))^{p+1} = (-\log(1-x))^{p} \big(\sum \frac{x^n}{n} \big)\\
 &\simeq \big(\sum p \frac{(\log n)^{p-1}}{n} x^n \big) \big(\sum
  \frac{x^n}{n}\big)\\
 & =\sum_n x^n ~\frac{p}{n} \Big( \sum_{k=1}^{n-1} \frac{(\log k)^{p-1}}{k} +
  \sum_{k=1}^{n-1} \frac{(\log k)^{p-1}}{n-k}\Big) \\
 &\simeq \sum_n x^n ~\frac{p}{n} \Big( \frac{(\log n)^{p}}{n} +(\log
  n)^p \Big) = (p+1)\sum_n x^n \frac{(\log n)^p}{n},   
\end{align*}
hence $S_{p,1,n} \sim \frac{1}{p!}(\log n)^p$.  The same inductive
proof works for the next-to-leading term of (\ref{eq:asymp-stirling}).
In that case, we also need the next-to-leading term of lemma
\ref{logkp-(n-k)}.  For each subsequent term that we want to prove in
(\ref{eq:asymp-stirling}), we need one more term of lemma
\ref{logkp-(n-k)} while the leading term of lemma \ref{logkp-k} is
enough.  One thus sees how the sequence of decreasing logarithms in
(\ref{eq:asymp-stirling}) is intimately related to that of $
\sum_{r=1}^{p} c_{p,r} (\log n)^{p-r} $ in lemma \ref{logkp-(n-k)}.

\section{Asymptotics of sums involving $1/k^q (\log k)^p$}
\label{sec:fourth-sum}

\begin{lem} \label{lem:1/logkp} For $p\geq 1$ we have:
  \begin{align*}
   \sum_{k=2}^{n-1} \frac{1}{(\log k)^p} =\tfrac{1}{(p-1)!} ~{\rm
  li}(n) -n ~\sum_{r=1}^{p-1} \frac{c_{p,r}}{(\log n)^r}
  +C_{p,0} - \frac{1}{2}\frac{1}{(\log n)^p} -\sum_{k=1}^m \Big(
 \sum_{r=1}^{2k-1} \frac{d_{p,r,k}}{(\log
  n)^{r+p}} \Big) \frac{1}{n^{2k-1}} +O(\tfrac
  {1}{n^{2m}}) ,
   \end{align*}
with $c_{p,r}:=\frac{(r-1)!}{(p-1)!}$ and $d_{p,r,k}
  := \frac{B_{2k}}{2k(2k-1)} \frac{(p-1+r)!}{(p-1)!} S_{r-1,1,2k-2}$,
  and $C_{p,0}$ is the constant term.\\
The log-integral is defined by li($z):= \int_0^z \frac{dt}{log t}$.
\end{lem}
\begin{proof}
  This follows from the Euler--Maclaurin formula with
  $f(x):=\frac{1}{(\log x)^p}$ and
{\scriptsize
  \begin{align*}
    \int f(x) &= \tfrac{1}{(p-1)!} ~{\rm li}(x) -x ~\sum_{r=1}^{p-1}
    \frac{(r-1)!}{(p-1)!} \frac{1}{(\log x)^r} \\
    f^{(i)}(x) &=(-1)^i (i-1)! \sum_{r=1}^i \frac{(p-1+r)!}{(p-1)!}
    \frac{S_{r-1,1,i-1}}{x^i~(\log x)^{r+p}}.
  \end{align*}
} From these, it is also straightforward to write down the `exact'
expression for the constant $C_{p,0}$, involving a formal (infinite)
sum over Bernoulli numbers. We omit it as it is not enlightening.
\end{proof}

Note that the second sum on the rhs is just the start of the
asymptotic expression of the first term, since li$(n)\approx n
\sum_{r\geq 1} \frac{(r-1)!}{(\log n)^r}$. So we might replace the two
terms by $ \frac{n}{(p-1)!} \sum_{r\geq p} \frac{(r-1)!}{(log n)^r}$.
This is indeed what one obtains when numerically looking for the
asymptotics of the lhs; the first term is $\frac{n}{(\log n)^p}$, and
correctly so.  Yet since this asymptotic expansion diverges for all
values $n$, the replacement would be disastrous for numerical
evaluation of the constant $C_{p,0}$.

\begin{lem}  \label{lem:1/klogkp} For $p\geq 1$ we have:
  \begin{align*}
\sum_{k=2}^{n-1} \frac{1}{k~ (\log k)^p} =C_{p,1}-\frac{1}{(p-1)(\log n)^{p-1}}
   -\frac{1}{2}\frac{1}{n(\log n)^p} -\sum_{k=1}^m \Big(
  \sum_{r=0}^{2k-1} \frac{d_{p,r,k}}{(\log
  n)^{r+p}} \Big) \frac{1}{n^{2k}} +O(\tfrac {1}{n^{2m+1}}),
  \end{align*}
with $d_{p,r,k} := \frac{B_{2k}}{2k} \frac{(p-1+r)!}{(p-1)!}
  S_{r,1,2k-1}$,  and $C_{p,1}$ is the 
  constant term. \\ For $p=1$, the second term on the rhs has to be
  replaced by $\log(\log n)$ (which becomes the leading term).
\end{lem}
\begin{proof}
  This follows from the Euler--Maclaurin formula with
  $f(x):=\frac{1}{x~(\log x)^p}$ and
{\scriptsize
  \begin{align*}
    \int f(x) &= -\frac{1}{(p-1)(\log x)^{p-1}} ~~~~~(p\geq 2) \\
    f^{(i)}(x) &=(-1)^i i! \sum_{r=0}^i \frac{(p-1+r)!}{(p-1)!}
    \frac{S_{r,1,i}}{x^{i+1}~(\log x)^{r+p}}.
  \end{align*}
}
\end{proof}

\begin{lem} For $p\geq 1$ and $q\geq 2$ we have:
  \begin{align*}
\sum_{k=2}^{n-1} \frac{1}{k^q~ (\log k)^p} &=C_{p,q}+
\tfrac{(1-q)^{p-1}}{(p-1)!}  {\rm Ei}({(1-q)\log n}) -\frac{1}{n^{q-1}}
\sum_{r=1}^{p-1} \frac{c_{p,q,r}}{(\log n)^r}
   -\frac{1}{2}\frac{1}{n^q (\log n)^p} \\
&\hspace{4cm} -\sum_{k=1}^m \Big(
  \sum_{r=0}^{2k-1} \frac{d_{p,q,r,k}}{(\log
  n)^{r+p}} \Big) \frac{1}{n^{2k-1+q}} +O(\tfrac {1}{n^{2m+q}}),
  \end{align*}
with $c_{p,q,r}:= \frac{(1-q)^{p-1-r} (r-1)!}{(p-1)!}$ and
 $d_{p,r,k} := \frac{B_{2k} (2k+q-2)!}{(2k)!(q-1)!} \frac{(p-1+r)!}{(p-1)!}
  S_{r,q,2k+q-2}$,  and $C_{p,q}$ is the constant term. 
The exponential integral function is defined by the principle value of
 the integral: {\rm Ei}$(x) := -\int_{-x}^\infty \frac{e^{-t}}{t} dt$.
\end{lem}
\begin{proof}
  This follows from the Euler--Maclaurin formula with
  $f(x):=\frac{1}{x^q~(\log x)^p}$ and
{\scriptsize
  \begin{align*}
    \int f(x) &= \frac{(1-q)^{p-1}}{(p-1)!}  {\rm Ei}((1-q)\log x)
    -\frac{1}{x^{q-1}} \sum_{r=1}^{p-1} \frac{
    (r-1)!}{(p-1)!} \frac{(1-q)^{p-1-r}}{(\log x)^r}\\
    f^{(i)}(x) &=(-1)^i \frac{(i+q-1)!}{(q-1)!} \sum_{r=0}^i
    \frac{(p-1+r)!}{(p-1)!} \frac{S_{r,q,i+q-1}}{x^{i+q}~(\log x)^{r+p}}.
  \end{align*}
}
\end{proof}
Again, the third term on the rhs is just the start of the asymptotic
expansion of the second term, since Ei$((1-q)n) \approx \sum_{r\geq 1}
\frac{(r-1)!}{(1-q)^r (\log n)^r} $.  So we might replace both terms
by the infinite sum $ n^{q-1} \sum_{r\geq p} \frac{c_{p,q,r}}{(\log
n)^r} $.  But since this diverges for all $n$, the replacement is
disastrous for numerically computing the constant $C_{p,q}$.

Note that when $p\geq 2$, the previous lemma makes sense also for
$q=1$, and one recovers the preceding lemma (since $(1-q)^{p-1-r}$
vanishes unless $r=p-1$).

For large $p$ or large $q$, it is quite obvious that the main
contribution to the sum $\sum_{k=2}^{n-1} \frac{1}{k^q~ (\log k)^p}$
comes from the term $k=2$ and that the constants $C_{p,q}$ will
converge towards $\frac{1}{2^q~ (\log 2)^p}$.  Just how quick they
converge can be empirically determined:
asymptotically for large $p$ or $q$, we have
$$
C_{p,q} \sim \frac{1}{2^q~ (\log 2)^p} + e^{-ap-bq} +e^{-cp-dq}+\dots,
$$
with $a=0.09405, b=1.0986, c=0.3266, d=1.386$.  Of course, these
values are nothing but $\log(\log 3)$, $\log 3$, $\log(\log 4)$, $\log
4$, so as to obtain $\frac{1}{3^q (\log 3)^p} +\frac{1}{4^q (\log 4)^p}$~!
So we come back from where we started. This comes as no surprise when
$C_{p,q}=\sum_{k=2}^{\infty} \frac{1}{k^q~ (\log k)^p}$, but it is a
surprise when the infinite sum does not converge, ie. when $C_{p,q}$
is not the leading term in the asymptotics, eg. when $q=0$ and 
$p$ becomes large.  

\begin{table*}[htbp]
  \centering
\hspace{.8cm}  \begin{tabular}[t]{r|llllll}
 $C_{p,q}$  & $q=0$ & 1 & 2 & 3 & 4  \\
\hline 
$p=1$:&      --0.24324 & 0.794679 & 0.605522 & 0.237996 & 0.106201 \\
 2:&      3.10329 & 2.10974 & 0.692606 & 0.305808 & 0.143463 \\
 3:&        4.96079 & 2.06589 & 0.882388 & 0.412914 & 0.199091 \\
 4:&       6.00344 & 2.55912 & 1.18928 & 0.573295 & 0.28066 \\
 5:&       7.46574 & 3.42982 & 1.65131 & 0.808652 & 0.399314 \\
 6:&        9.92015 & 4.75831 & 2.33023 & 1.15106 & 0.571244 
  \end{tabular}
  \begin{tabular}[t]{r|llllll}
 $2^{-q} (\log 2)^{-p}$  & $q=0$ & 1 & 2 & 3 & 4  \\
\hline 
$p=1$:&       1.4427 & 0.721348 & 0.360674 & 0.180337 & 0.0901684 \\
 2:&         2.08137 & 1.04068 & 0.520342 & 0.260171 & 0.130086 \\
 3:&         3.00278 & 1.50139 & 0.750695 & 0.375348 & 0.187674 \\
 4:&         4.3321 & 2.16605 & 1.08302 & 0.541512 & 0.270756 \\
 5:&         6.24989 & 3.12495 & 1.56247 & 0.781237 & 0.390618 \\
 6:&         9.01669 & 4.50835 & 2.25417 & 1.12709 & 0.563543 
  \end{tabular}
  \caption{Comparison between $C_{p,q}$ and $\frac{1}{2^q (\log 2)^p}$.}
   \label{tab:table1}
\end{table*}
We may want to add li$(2)=1.045$ to $C_{1,0}$  so as to obtain the constant\\
lim$_{n\to \infty} \Big( \sum_{k=2}^{n-1} \frac{1}{\log k} -\int_2^n
\frac{dx}{\log x}\Big) =0.80192543$. Similarly, we add $\log(\log
2)=-0.36651$ to $C_{1,1}$, so as to obtain the constant
lim$_{n\to \infty} \Big( \sum_{k=2}^{n-1} \frac{1}{k(\log k)}
-\int_2^n \frac{dx}{x(\log x)}\Big) =0.4281657$.  Both values already
occurred in \cite{B-77}, see also \cite{F-03}.  We were not able to
recognize an exact form for either of these two constants (using
PARI for integer linear combinations of other constants, or using
Plouffe's inverter or his Maple code).

\section{Asymptotics of the Taylor coefficients of $(-z/\log(1-z))^k$}

We now use lemma \ref{lem:1/klogkp} to generalise a result known to
P\'olya \cite{P-54} about the Taylor coefficients of a certain
generating function. In 1954, P\'olya \cite{P-54} noted that
\begin{equation} \label{eq:polya}
a_n \sim -\frac{1}{n (\log n)^2} \qquad\textrm{for} \qquad
f(z)=\frac{z}{-\log(1-z)} =: \sum a_n z^n.  
\end{equation}
We shall be interested in the asymptotics of the $a_n$ when the
generating function is raised to some power $k$ (positive integer).
For $k=1$, the series begins as $1-\half x -\tfrac{1}{12} x-\dots$ and
all coefficients are negative except $a_0$.  The $a_n$ for $k=2$ are
asymptotically given by the convolution of those at $k=1$, viz.  $a_n
=\sum \frac{1}{i(\log i)^2 ~(n-i)(\log (n-i))^2 }$. Since this sum
makes only sense for $i$ running from 2 to $n-2$, we write the terms
$-\frac{1}{n(\log n)^2} - \frac{-1/2}{(n-1)(\log (n-1))^2} \simeq
-\frac{1/2}{n(\log n)^2}$ twice separately. By symmetry, we can write:
$$
a_n \approx -\frac{1}{n(\log n)^2} + 2\sum_{i=2}^{(n-2)/2} \frac{1}{i(\log
  i)^2 ~(n-i)(\log (n-i))^2 } ,
$$
where in the last sum, $\log(n-i) \geq \log n/2$ and
$\frac{1}{i(n-i)}= \frac{1/n}{i}+\frac{1/n}{n-i}$. From lemma
\ref{lem:1/klogkp} we know that $\sum^{n/2} \frac{1}{i(\log i)^2} =C
-\frac{1}{(\log n/2)} +O(\frac{1}{(\log n/2)^2})$, wherein the constant
$C$ is figurative, since the quantities  $1/(i(\log i)^2)$ only
approximate the exact values of the Taylor coefficients.  Further,
$\sum^{n/2} \frac{1}{(n-i)(\log i)^2} \leq \frac{1}{n/2} \sum^{n/2}
\frac{1}{(\log i)^2} =O(\frac{1}{(\log n)^2})$.  Overall:
$$
a_n \approx -\frac{1}{n(\log n)^2} +\frac{2C}{n(\log n)^2}
-\frac{2}{n(\log n)^3}  +O\Big(\frac{1}{n(\log n)^4}\Big) 
$$
Since the singularity of $f(z)^2$ at $z=1$ is of higher order than
that of $f(z)$, the decrease of coefficients should be stronger; hence
the $\frac{1}{n(\log n)^2}$ terms have to cancel each other and so
$C=1/2$. We are left with $a_n \approx -\frac{2}{n(\log
  n)^3}+\dots$. One similarly obtains:
\begin{equation}\label{eq:polya2}
a_n \approx -\frac{k}{n(\log n)^{k+1}} +O\Big(\frac{1}{n(\log n)^{k+2}}\Big)
\qquad\textrm{for} \qquad f(z)=\Big(\frac{z}{-\log(1-z)}\Big)^k.  
\end{equation}

A naive attempt at justifying P\'olya's result (\ref{eq:polya}) would be to
use Cauchy's formula $a_n = \frac{1}{2\pi i} \oint_C \frac{f(z)
  ~dz}{z^{n+1}}$ and to compute the contour integral on the unit
circle, $z=e^{i\theta}$. Note that $-\log(1-e^{i\theta}) =-\log(-2i
e^{i\theta/2} \sin\frac{\theta}{2}) $. Thus we would have (wrongly)
$$
a_n = \frac{1}{2\pi} \int_0^{2\pi}
\frac{e^{-i(n-1)\theta} d\theta}{i\frac{\pi-\theta}{2}
  -\log(2\sin\frac{\theta}{2})} 
\simeq  \frac{1}{2\pi}  \int_0^{c/n} \frac{d\theta}{-\log\theta}
~~~~~\approx \frac{c/(2\pi)}{n \log n} +O\Big( \frac{1}{n (\log n)^2}\Big)
$$
where we replaced $ \int_0^{2\pi}$ by $ \int_0^{c/n}$ since for
large $n$, only small values of $\theta$ will contribute substantially
to the integral.  This approximation, however, does not yield the
desired result -- presumably because one cannot replace
$e^{-i(n-1)\theta}$ by 1.  Similarly, had we used partial integration
with $f'(\theta)= \frac{e^{i\theta} \frac{1}{2}\cot\frac{\theta}{2}}
{(i\frac{\pi-\theta}{2} -\log(2\sin\frac{\theta}{2}))^2}$, we would
have ended up with $ \frac{1}{2\pi i n} \int_0^{c/n}
\frac{e^{-in\theta} d\theta}{\theta (\log\theta)^2} \approx
\frac{c/(2\pi)}{n \log n}$, again with the wrong leading term.  As the
integral is not tractable by the Laplace method, the saddle point
method or any other trick described in \cite{dB-58}, we shall see in
the next subsection that the solution lies in a clever choice of the
contour of integration.

The result (\ref{eq:polya2}) is not new, but was already obtained by
N\"orlund in 1961 using combinatorics of Bernoulli polynomials, and
rederived by Flajolet and Odlyzko in 1990 by evaluating the contour
integral in Cauchy's formula.  For completeness, we present these two
alternative and elegant paths below.

\subsection{The Flajolet--Odlyzko approach}

In 1990, Flajolet and Odlyzko summerised the `transfer properties' of
analytic functions, viz. the behaviour of the function at the first
singularity on the convergence radius is directly reflected in the
behaviour of the Taylor coefficients.  One of their result is \cite{FO-90}:

\begin{thm} (Flajolet--Odlyzko, 1990)
Let $f(z)$ be analytic in $|z| < 1+\eta$ except for a singularity at
$z=1$, and let 
$$
f(z) = O\Big( (1-z)^\alpha (-\log(1-z))^\gamma \Big) ~~~\textrm{as}~~~
z\to 1 ~~~(\alpha, \gamma \in\R).
$$
Then the coefficients $a_n$ in $f(z)=\sum a_n z^n$ grow like
$a_n=O\Big(\frac{(\log n)^\gamma}{n^{\alpha+1}}\Big) $.  
\end{thm}
\begin{proof} (sketchy).
It is comforting to see that the proof boils down to a mere
application of the Cauchy formula, ie. a contour integral around the
origin, viz. $a_n= \frac{1}{2\pi i} \oint_C \frac{f(z)
  ~dz}{z^{n+1}}$, but the contour has to be chosen cleverly -- as in
figure~\ref{fig:pic1}.
\begin{figure}[htbp]
  \begin{center}
    \begin{picture}(0,75)  
\setlength{\unitlength}{.8pt} 
\put(0,50){
\arc[100](50,10){337}  
\put(35,0){\arc[100](5,5){273}}  
\curve[100](40,5, 50,10)
\curve[100](40,-5, 50,-10)
\put(33,-3){1}
\put(0,-2){0}
\put(23,12){$C_1$}
\put(-45,7){$C$}
}
    \end{picture}
  \end{center}
\caption{The contour of integration, excluding the singularity at $z=1$.}
  \label{fig:pic1}
\end{figure}
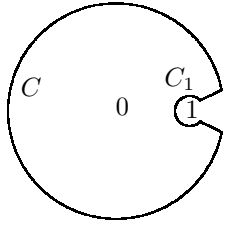

The contour $C$ will be a circle of radius $1+\eta$ with a tiny
roundabout around the singularity; the main contribution will come
from this little near-circle $C_1$ (of radius $1/n$) around $z=1$.
Note that for any compact domain inside our contour $C$, there is a
constant $K$ such that $ |f(z)|\leq K | (1-z)^\alpha
(-\log(1-z))^\gamma | $. On the circle $C_1$, we have
$1-z=e^{i\theta}/n$ and the following bounds: $|f|\leq
K(\frac{1}{n})^\alpha $ sup$|\log(n e^{-i\theta})|^\gamma = O\Big(
\frac{(\log n)^\gamma}{n^{\alpha}} \Big)$, as well as $|z|^{n+1}\geq
(1-\frac{1}{n})^{n+1} \to 1$ and $\int_{C_1} |dz| \leq 2\pi/n$. Hence
the main contribution to the contour integral can be estimated by:
$a_n \approx \frac{1}{2\pi} \oint_{C_1} \frac{|f(z)| ~|dz|}{|z|^{n+1}}
= \Big(\frac{(\log n)^\gamma}{n^{\alpha+1}}\Big) $, as we wished.
\end{proof}

This result was readily obtained, but is treacherous when $\alpha$ is
a non-negative integer, say 0: the Taylor coefficients of $-\log(1-z)$
decrease like $1/n$ and not like $(\log n)/n$, and those of
$1/(-\log(1-z))$ decrease like $1/(n (\log n)^2)$ and not like $1/(n
(\log n))$. For this case, it is useful to have a precise asymptotic
development, whose derivation we sketch as follows (see \cite{FO-90}
for details).  Let $f(z)$ be the function $(1-z)^\alpha (-\frac{1}{z}
\log (1-z))^\gamma$.  Change to the variable $z=1+t/n$ and expand
$\log(-n/t)^\gamma = (\log n)^\gamma \sum_{i\geq 0} \binom{\gamma}{i}
\big( -\frac{\log (-t)}{\log n}\big)^i$.  The contour integral now
contains the piece $G_i :=\frac{1}{2\pi i} \int_{C_2} (-t)^\alpha
(\log(-t))^i e^{-t} dt$, whose contour can be deformed to a well-known
integral that simply yields $G_i= \p_\alpha^i
\frac{1}{\Gamma(-\alpha)}$.  Hence the main contribution from the
Cauchy formula for $f(z)=\sum a_n z^n$ is:
\begin{equation}\label{eq:main-asymp}
a_n \approx \frac{(\log n)^\gamma}{n^{\alpha+1}} \sum_{i\geq 0} \frac{(-1)^i
  \binom{\gamma}{i} G_i}{(\log n)^i} \qquad \textrm{for}\qquad f(z)=
(1-z)^\alpha (-\frac{1}{z} \log (1-z))^\gamma.
\end{equation}
The first term of the sum, $i=0$, yields $\frac{1}{\Gamma(-\alpha)}$
and thus simply drops out when $\alpha$ is a non-negative integer
(where $\Gamma(-\alpha)$ has a pole).  This explains the above treachery.
For $\alpha=0$ and $\gamma=-k$ (negative integer), we have
$\binom{-k}{i} =(-1)^i \binom{k+i-1}{i}$ and $g_i :=\p_\alpha^i \big|_0
\frac{1}{\Gamma(-\alpha)} =-1,2\gamma,-3\gamma^2+\frac{\pi^2}{2},\dots
(i=1,2,\dots)$, and so obtain our sought-for asymptotics:
\begin{equation}\label{eq:main-asymp2}
a_n \approx \sum_{i\geq 1} \frac{ \binom{k+i-1}{i} g_i}{n~(\log
  n)^{i+1}} = \frac{-k}{n(\log n)^2} +\dots \qquad\qquad \textrm{for}\qquad
f(z)= \Big( \frac{z}{-\log (1-z)} \Big)^k.  
\end{equation}

\subsection{N\"orlund's approach}

Most surprisingly, (\ref{eq:main-asymp}) was arrived at 20 years
earlier from quite a different angle, namely by N\"orlund \cite{N-61}
who recognized the Taylor coefficients $a_n$ of $f(z)$ as being values
of generalised Bernoulli polynomials:
\begin{equation} \label{eq:bernoulli-poly}
(1-z)^\alpha (-\frac{1}{z} \log (1-z))^\gamma = \sum_{n\geq 0}
\frac{(-z)^n}{n!} B_n^{(n+\gamma+1)}(\alpha).  
\end{equation}
The Bernoulli polynomials of order $\gamma$ are defined by
$$
\Big( \frac{t}{e^t -1}\Big)^\gamma e^{\alpha t} =: \sum_{n\geq 0}
B_n^{(\gamma)}(\alpha) \frac{t^n}{n!}
$$
and coincide with the usual Bernoulli polynomials for $\alpha=0$
(and thus yield the Bernoulli numbers for $\gamma=1$).  To arrive at
(\ref{eq:bernoulli-poly}), use Cauchy's formula $
\frac{B_n^{(\gamma)}(\alpha)}{n!} = \frac{1}{2\pi i} \oint
\frac{t^{\gamma-n-1} e^{\alpha t}}{(e^t -1)^\gamma} dt$, then
substitute $t=\log(1-z)$ and shift $\gamma\to \gamma+n+1$; you thus
obtain the Cauchy formula equivalent to (\ref{eq:bernoulli-poly}).
At $\alpha=1$ and $\gamma\geq0$, we recover the Stirling numbers of
first kind: 
$$
(-1)^n ~B_n^{(n+\gamma+1)}(1) = \Big[ {n+\gamma \atop \gamma}
  \Big] \Big/ \binom{n+\gamma}{\gamma} = \frac{n!\gamma!}{n+\gamma}
~S_{\gamma-1,1,n+\gamma-1} 
$$
The polynomials satisfy 
\begin{align*}
B_n^{(\gamma)}(\alpha) =\int_\alpha^{\alpha+1} B_n^{(\gamma+1)}(t) ~dt 
\qquad \textrm{and} \qquad B_n^{(n+1)}(\alpha)
&=(\alpha-1)(\alpha-2)\cdots(\alpha-n)\\
& = (-1)^n (n-\alpha)\cdots
(1-\alpha) =  (-1)^n \frac{\Gamma(n-\alpha+1)}{\Gamma(1-\alpha)}
\end{align*}
Hence we can readily obtain an asymptotic expression for \underline{$\alpha=0$
and $\gamma=-1$}:
\begin{align*}
\frac{(-1)^n ~B_n^{(n)}(\alpha)}{n!} &= \int_0^1
\frac{\Gamma(n+1-\alpha-t)}{\Gamma(1-\alpha-t) \Gamma(n+1)} ~dt \\
&=  \int_0^1 n^{-\alpha-t} \Big( \frac{1}{\Gamma(1-\alpha-t)} + O\big(
  \tfrac{1}{n}\big) \Big) \\
&= n^{-\alpha} \Big[ \frac{1}{(\log n) \Gamma(1-\alpha)} -\int_0^1
\frac{e^{-t \log n}}{-\log n} ~\p_t  \frac{1}{\Gamma(1-\alpha-t)} dt +
O\big( \tfrac{1}{n \log n}\big) \Big] 
\end{align*}
where we used Stirling's approximation $n! \sim e^{n\log n
  -n+\half\log (2\pi n)}$ in the second step and partial integration
in the third step.  Further partial integrations yield
$$
\frac{(-1)^n ~B_n^{(n)}(\alpha) ~n^\alpha}{n!} = \sum_{i=0}^{r-1}
\frac{1}{(\log n)^{i+1}} ~\p_\alpha^i \frac{1}{\Gamma(1-\alpha)} +O\big(
\tfrac{1}{(\log n)^{r+1}} \big)
$$
where the derivative term evaluated at $\alpha=1$ is:
$\p_\alpha^i\big|_1 \frac{1}{\Gamma(1-\alpha)} =
\p_\alpha^i\big|_0 \frac{1}{\Gamma(-\alpha)} =:g_i$, with
$g_i=0,-1,2\gamma,-3\gamma^2+\frac{\pi^2}{2},\dots$ for
$i=0,1,\dots$~.  Overall: 
$$
a_n =\frac{(-1)^n ~B_n^{(n)}(1)}{n!} \approx \sum_{i\geq 1}
\frac{ g_i}{n (\log n)^{i+1}} \qquad\textrm{for} \qquad f(z)=
\frac{z}{-\log (1-z)} . 
$$
From here, it is straightforward (via induction) to deduce the general
case $\gamma=-k$ (negative integer):
$$
a_n =\frac{(-1)^n ~B_n^{(n-k+1)}(1)}{n!} \approx \sum_{i\geq 1}
\frac{\binom{k+i-1}{i} g_i}{n (\log n)^{i+1}} = \frac{-k}{n(\log n)^2} +\dots
\qquad\textrm{for} \qquad f(z)= \Big( \frac{z}{-\log (1-z)} \Big)^k
$$ 
in full agreement with (\ref{eq:main-asymp2})

\section{Asymptotics of the incomplete Gamma function}

This sections examines whether the constant $C_{1,1}$ met in section 7
also occurs in the constant term of the the function
$\sum_{n=2}^\infty \frac{1}{(n\log n)^s}$.  The answer is yes, as the
general theory shows, but $\gamma$ has to be subtracted to obtain the
full constant term.  The resulting expansion is presented in
(\ref{eq:meta-zeta}), and on the way we shall derive the following
intermediate result about the asymptotics of the incomplete gamma
function ($a$ is an arbitrary constant):
$$
\Gamma(-s,as) := \int_{as}^\infty t^{-s-1} e^{-t} dt = -\log (as) -\gamma
+s [\half (\log (as))^2 + a -c_1 ] + O(s^2 \log^3(as)) \qquad\textrm{as
}~s\to 0.
$$

\subsection{Preliminaries}

Recall the coincidence in the constant terms of the following
asymptotic expansions: 
\begin{align*}
  \zeta(s) =\sum_{k\geq 1} \frac{1}{k^s} &\approx \frac{1}{s-1} +\gamma
  -\gamma_1 (s-1)+\dots \qquad (s\to 1)\\ 
  \sum_{k=1}^{n-1} \frac{1}{k} &\approx \log n + \gamma -\frac{1}{2n}+\dots
 \qquad (n\to \infty)
\end{align*}
Landau confirmed this coincidence for a broader class of Dirichlet
series: suppose $ \sum_{n\leq x} h(n) \sim \alpha x+\dots$ (among
other constraints on $h(n)$), then:
$$
 \sum_{n=1}^\infty \frac{h(n)}{n^s} \approx \frac{\alpha}{s-1} +\beta+\dots
 \qquad (s\to 1^+)
$$
for some constant $\beta$, and
$$
\sum_{n\leq x} \frac{h(n)}{n} \approx \alpha \log x + \beta +\dots \qquad
(x\to \infty). 
$$

We shall be concerned with a weaker generalisation.  First recall the
discrete partial integration formula for some continuous function
$\phi$ and some sequence $a_n$ with primitive $A(x):=\sum_{n\leq x} a_n$:
$$
\sum_{n=a+1}^b a_n \phi(n) = A(x) \phi(x) \Big|_a^b -\int_a^b A(x)\phi'(x)dx.
$$
When $a,b \in\Z$ and $a_n=1$ with $A(x)=\lfloor x \rfloor$, the
formula reduces to
\begin{align*}
\sum_{n=a}^b \phi(n) &=x\phi(x)\Big|_a^b -\int_a^b \lfloor x \rfloor \phi'(x)dx
+\phi(a)\\
&= \int_a^b \phi + \int_a^b (x- \lfloor x \rfloor) \phi' + \phi(a) .
\end{align*}
If $\phi_s(x)$ is a suitable function depending on a parameter $s$,
like $\phi_s(x)=\frac{1}{x^s}$, the derivative in the second integrand
will ensure that we may exchange the limits $b\to\infty$ and $s\to 1$
(since $x- \lfloor x \rfloor$ is bounded).  In other words, the
first integral contains the singularity as $s\to 1$, while the second
integral yields merely a constant.  Denote by $\phi_1$ the
function $\phi_s$ obtained after taking the limit $s\to1$; we then
perform the partial integration backwards:
$$
\int_a^b (x- \lfloor x \rfloor) \phi_1'  = \sum_{k=a}^{b-1}
\int_k^{k+1}(x-k) \phi_1'  = \sum_{k=a}^{b-1} \Big(\phi_1(k+1)
-\int_k^{k+1} \phi_1 \Big) =  \sum_{k=a+1}^{b} \phi_1(k) -\int_a^b \phi_1.
$$
Having previously taken the limit $b\to\infty$, we obtain our desired
generalisation: 
$$
\lim_{s\to 1} \Big( \sum_{n=a}^\infty \phi_s(n) -\int_a^\infty \phi_s \Big) 
= \lim_{b\to \infty}  \Big( \sum_{k=a}^{b} \phi_1(k) -\int_a^b \phi_1 \Big).
$$
For instance, for $\phi_s(x)=\frac{1}{x^s}$ and $a=1$, we have:
$$
\lim_{s\to 1} \Big( \zeta(s) -\frac{1}{s-1} \Big) 
= \lim_{b\to \infty}  \Big( H_b -\log b \Big) = \gamma.
$$
For $\phi(x)=\frac{1}{(x\log x)^s}$ and $a=2$, we have  by lemma
\ref{lem:1/klogkp}:
\begin{equation}\label{eq:1/nlogns}
\lim_{s\to 1} \Big( \sum_{n=2}^\infty \frac{1}{(n\log n)^s}
-\int_2^\infty \frac{dx}{(x\log x)^s} \Big)  
= \lim_{b\to \infty} \Big( \sum_{n=2}^{b} \frac{1}{n\log n} -\log(\log
b)+\log(\log 2) \Big) = C_{1,1} +\log(\log 2).
\end{equation}

\subsection{Gamma function asymptotics}

In order to find the constant term in the asymptotics of the function
$ \sum_{n=2}^\infty \frac{1}{(n\log n)^s}$, we still need to expand
the corresponding integral up to the constant term.  This will include
finding its singularity at $s=1$. Note first that
$$
\int_2^\infty \frac{dx}{(x\log x)^s} = (s-1)^{s-1}
~~\Gamma(1-s,~(s-1)\log 2) = (s-1)^{s-1} \int_{(s-1)\log 2}^\infty
t^{-s} e^{-t} dt.
$$
The incomplete Gamma function is defined by $\Gamma(z,s)
:=\int_s^\infty t^{z-1}e^{-t} dt$.  Hence $\frac{d}{ds} \Gamma =
\Gamma_z \frac{dz}{ds} +\Gamma_s $.  To simplify matters, we shall
first study the behaviour of $\Gamma(-s,s)$ as $s\to 0$.  Note that
$\Gamma(0,s) \sim -\log s -\gamma+\dots $ and $\Gamma(-s,0) \sim
-\frac{1}{s}-\gamma+\dots$, so that the simultaneous limit will behave
like the weaker of both, ie lim$_{s\to 0} \Gamma(-s,s) \sim -\log s +\dots $.
Equipped with this intuition, we proceed by expanding $\p_s
\Gamma(-s,s)$ (around $s=0$) and then integrating the expanded result. 
Now with $z(s)=-s$ we have
$$
\frac{d}{ds} \Gamma(-s,s) = \bigg[ \frac{\p\Gamma}{\p z} \frac{dz}{ds}
+\frac{\p\Gamma}{\p s} \bigg]_{z=-s} =-e^{-s} s^{-s-1} - \int_s^\infty
t^{-s-1} e^{-t} (\log t) dt.
$$
The first part is expanded as $-e^{-s} s^{-s-1} = -\frac{1}{s}
e^{-s(1+\log s)} = -\frac{1}{s} + (\log s +1) +\dots$ and gives us --
upon integration  -- the required $(-\log s)$ term, so that the second
part can contain at most log-singularities.
Differentiating the latter gives us $ e^{-s} s^{-s-1} \log s
+\int_s^\infty  t^{-s-1} e^{-t} (\log t)^2 dt $, which behaves as
$(\log s)/s +\dots$, so that the original integral behaves as $\half
(\log s)^2 +$ const $+\dots$.  In general, $\int_s^\infty  t^{-s-1}
e^{-t} (\log t)^i dt $ will behave as $-\frac{1}{i+1}(\log s)^i +$
const $+\dots$, and by bootstrapping we obtain: 
\begin{align*}
\int_s^\infty  t^{-s-1} e^{-t} dt  &= -(\log s) +c_0 + 
  s [\half(\log s)^2 + 1-c_1 ] +\dots \\
\int_s^\infty  t^{-s-1} e^{-t} (\log t) dt &= -\half(\log s)^2 +c_1 +
  s [\tfrac{1}{3}(\log s)^3 +(\log s) -1-c_2 ] +\dots \\
\vdots\\
\int_s^\infty  t^{-s-1} e^{-t} (\log t)^i dt &= -\tfrac{1}{i+1}(\log
  s)^{i+1} +c_i + s \Big[ \tfrac{1}{i+2}(\log s)^{i+2} +\sum_{r=0}^i
   \tfrac{(-1)^{r+i} i!}{r!}(\log s)^{r} -c_{i+1} \Big] +\dots
\end{align*}
In this way, we can start with the $n$-th line at order $O(s)$ and
recursively determine the expansion of the first line up to $O(s^n)$
in terms of the constants of integration $c_0,\dots,c_n$.  The latter
can be empirically determined: $c_0=-\gamma$, $c_1=0.98905$, $c_2=
-1.81497$, $c_3=5.89038$, $c_4=-23.568$, etc. We have computed the
first 300 of them and it seems that their asymptotics are $ c_i =
(-1)^{i+1} (i! -e^{0.8~i\log i} +\dots) $. It would be interesting to
know more about these constants. In particular, we obtain our desired
expansion:
$$
\Gamma(-s,s) =  -(\log s) -\gamma + s [\half(\log s)^2 + 1-c_1 ] +
  s^2 \Big[ -\tfrac 16 (\log s)^3 -(\log s) +\tfrac 34 +
  \tfrac{c_2}{2} \Big] +\dots \qquad (s\to 0)
$$
Next, study $\Gamma(-s,a s)$ for some constant $a$.
Now $\frac{d}{ds} \Gamma = \Gamma_z (-1)+\Gamma_s a$. 
Bootstrapping yields now:
$$
\int_{as}^\infty  t^{-s-1} e^{-t} (\log t)^i dt = -\tfrac{1}{i+1}(\log
  s)^{i+1} +c_i + s \Big[ \tfrac{1}{i+2}(\log s)^{i+2} +a \sum_{r=0}^i
   \tfrac{(-1)^{r+i} i!}{r!}(\log s)^{r} -c_{i+1} \Big] +\dots,
$$
with the constants $c_i$ taking the same values as before. In particular,
$$
\Gamma(-s,as) =  -\log (as) -\gamma +s [\half (\log (as))^2 + a -c_1 ] +
  s^2 \Big[ -\tfrac 16 (\log (as))^3 -a\log(as) +
  \tfrac{c_2}{2} -\tfrac{a^2}{4}+a \Big] +\dots\quad (s\to 0)
$$
Thus we can answer the question above:
\begin{align*}
\int_2^\infty \frac{dx}{(x\log x)^s} &= \Big[ 1+
(s-1)\log(s-1)+\dots\Big] \Big[ -\log ((s-1)\log 2) -\gamma
+(s-1)[\half\log^2((s-1)\log2)+\log2-c_1]+\dots \Big] \\
&=-\log(s-1) -(\gamma+\log\log 2) +(s-1)\Big[ -\half\log^2(s-1)
+(\log2-c_1-\gamma) \log(s-1) +\\
&\hspace{7cm} +\log2-c_1 \Big] + O((s-1)^2\log^3(s-1))
\end{align*}
Recalling (\ref{eq:1/nlogns}), we deduce the constant term in the
asymptotic development of the original function:
\begin{equation}
  \label{eq:meta-zeta}
\sum_{n=2}^\infty \frac{1}{(n\log n)^s} \approx -\log(s-1) +C_{1,1}-\gamma
 +O\big((s-1)\log^2(s-1)\big) \qquad (s\to 1)  
\end{equation}
where $C_{1,1}=0.794679...$ was given in table \ref{tab:table1}.

\section{A representation of polylogs and of Nielsen integrals}
\label{sec:polylogs}

In this section we present a representation of polylogs in terms of
Nielsen integrals which we have come across while embarking onto the
proofs in section~6.  Conversely, in (\ref{eq:S11}) and (\ref{eq:S12})
we give representations of Nielsen integrals involving Bernoulli
numbers, which boil down to mysterious identities for Bernoulli
numbers and harmonic numbers (conjecture~\ref{conj:bernoulli}).  We
only treat the cases $ \mathfrak{S}_{1,1}(x)$ and $
\mathfrak{S}_{1,2}(x)$, but are convinced that similar formulae hold
for all $ \mathfrak{S}_{1,p}(x)$.

\begin{lem}
  \begin{align}
    \Li_2(1-x) &= \zeta(2) +(\log x) \sum_1 \frac{x^n}{n} -\sum_1
    \frac{x^n}{n^2} \\
    \Li_3(1-x) &= \zeta(3) -\zeta(2) \sum_1 \frac{x^n}{n} - (\log x)
    \sum_1 \frac{H_{n-1}}{n} x^n + \sum_1 \Big(
    \frac{H^{(2)}_{n-1}}{n} + \frac{H_{n-1}}{n^2} \Big) x^n\\
    \Li_4(1-x) &= \zeta(4) -\zeta(3) \sum_1 \frac{x^n}{n} +\zeta(2)
    \sum_1 \frac{H_{n-1}}{n} x^n + (\log x)\sum_1
    \frac{\sum_{i=1}^{n-1} \frac{H_{i-1}}{i}}{n} x^n - \nonumber \\
    & \hspace{3cm} -\sum_1 \Big(
    \frac{\sum_{i=1}^{n-1} (\frac{H^{(2)}_{i-1}}{i}
    +\frac{H_{i-1}}{i^2} )}{n} +\frac{\sum_{i=1}^{n-1}
    \frac{H_{i-1}}{i} }{n^2} \Big) x^n \label{eq:Li_4(1-x)} \\
    \Li_j(1-x) &= \zeta(j) -\zeta(j-1) \sum_1 \frac{S_{0,1,n-1}}{n}
    x^n +\zeta(j-2) \sum_1 \frac{S_{1,1,n-1}}{n} x^n - \dots +(-1)^j
    \zeta(2) \sum_1 \frac{S_{j-3,1,n-1}}{n} x^n \nonumber\\ 
   & ~~+(-1)^j (\log x) \sum_1 \frac{S_{j-2,1,n-1}}{n} x^n -(-1)^j \sum_1
    \Big( \frac{T_{j-2,n-1}}{n} +\frac{S_{j-2,1,n-1}}{n^2} \Big) x^n
    \label{eq:Li_j(1-x)} 
  \end{align}
  with $T_{j-2,n-1} := \sum \frac{1}{i_1 \cdots i_{j-3} k^2}$ (sum
  over all $1 \leq i_1 <\cdots < i_{j-3} \leq n-1 $ and $1\leq k\leq
  n-1$, $k\neq i_1,\dots, i_{j-3}$) satisfying the following
  recursion: $ T_{j-2,n-1} = \sum_{i=1}^{n-1} \frac{T_{j-3,i-1}}{i} +
  \frac{S_{j-3,1,i-1}}{i^2} $.  The last sum of (\ref{eq:Li_j(1-x)})
  can also be written as $-\sum_{r=0}^{j-2} \mathfrak{S}_{1,j-r-1}(x)
  \frac{\log^{r}(1-x)}{r!}$ with $\mathfrak{S}_{1,p}(x) :=
  \sum_{n\geq 1} \frac{S_{p-1,1,n-1}}{n^2} x^n$, that is\footnote{We
    have used $(-\log(1-x))^p =p! \sum_{n\geq 0} \big[ {n \atop p}
    \big] \frac{x^n}{n!} = p! \sum_{n\geq 0} \frac{S_{p-1,1,n-1}}{n}
    x^n$ which is proved by expanding: lhs $= \p_r^p \big|_{_0}
    (1-x)^{-r} =\dots$ (notation and convention from \cite{GKP-89})} $x\p_x
  ~\mathfrak{S}_{1,p}(x) = \frac{(-\log(1-x))^p}{p!}$ for $p\geq 1$.
  Thus the above can be rewritten as ($j\geq 2$)
  \begin{equation} \label{eq:Li_j(1-x)bis}
    \Li_j(1-x) = \sum_{r=0}^{j-1} \Big( \zeta(j-r) -
    \mathfrak{S}_{1,j-r-1}(x) \Big) \frac{\log^{r}(1-x)}{r!}
  \end{equation}
wherein the term with $\zeta(1)$ should be dropped.
\end{lem}
\begin{proof}
  Show the recursion $\frac{d}{d \log(1-x)} \Li_j(1-x) =\Li_{j-1}(1-x)$
  using the fact that\\ $\sum_{r=0}^{j-1}  (\p_x
  \mathfrak{S}_{1,j-r-1}(x))  \frac{\log^{r}(1-x)}{r!} =
  \frac{\log^{j-1}(1-x)}{x} \sum_{r=0}^{j-1}
  \frac{(-1)^{j-1-r}}{(j-1-r)! r!} =0$ (binomial formula for $(1-1)^{j-1}$).
Note that $ \frac{\log^{r}(1-x)}{r!}$ could be replaced by
  $\frac{(-\Li_1(x))^r}{r!}$ or by $(-1)^r
  ~\Li_{1^r}(x)$, where the generalised polylog
  is defined by \\
$\Li_{s_1,\dots,s_k}(x) := \displaystyle \sum_{n_1>\dots >n_k>0}
  \frac{x^{n_1}}{n_1^{s_1} \dots n_k^{s_k}}$.  Thus the formula could
  read: 
$$ 
\Li_j(1-x) = \sum_{r=0}^{j-1} \big( \Li_{j-r}(1) -
    \Li_{2,1^{j-r-2}}(x) \big) (-1)^r ~\Li_{1^r}(x),
$$
wherein the term  with  $\Li_1(1)$ should be dropped.

\end{proof}
The $\mathfrak{S}_{1,p}$ are special cases of so-called {\em Nielsen
  integrals} (for $k\geq 1$): 
$$
\mathfrak{S}_{k,p}(x) := \frac{(-1)^{k+p-1}}{(k-1)! p!}
\int_0^1 \frac{dy}{y} \log^{k-1}(y) \log^p(1-xy) = \sum_{n\geq 1}
\frac{S_{p-1,1,n-1}}{n^{k+1}} x^n  = \Li_{k+1,1^{p-1}}(x) ,
$$
so that $\mathfrak{S}_{k-1,1}(x)
= \Li_k(x)$ and $ x \p_x ~\mathfrak{S}_{k,p}(x) =
\mathfrak{S}_{k-1,p}(x)$.  Also: $\mathfrak{S}_{1,1}(x)
=\Li_2(x)$, $\mathfrak{S}_{1,0}(x) =\log x$.

As a corollary, equate the last term of (\ref{eq:Li_4(1-x)}) with the
corresponding quantity in (\ref{eq:Li_j(1-x)bis}) and find: $ \frac{1}{n} 
\sum_{i=1}^{n-1} \Big( \frac{H_{i-1}^{(2)}}{i} +\frac{H_{i-1}}{i^2} \Big)
= \sum_{i=1}^{n-1}  \Big( \frac{H_{i-1}}{i (n-i)^2}
-\frac{H_{i-1}}{i^2 (n-i)} \Big)$; or, using partial fractions:
$$
\sum_{i=1}^{n-1} \frac{H_{i-1}}{(n-i)^2} = \sum_{i=1}^{n-1} \Big(
 \frac{H_{i-1}^{(2)}}{i} +2 \frac{H_{i-1}}{i^2} \Big) ,
$$
which can easily be generalised:
\begin{lem} For $p \geq 1$:
 $$
\sum_{i=1}^{n-1} \frac{H_{i-1}}{(n-i)^p} = \sum_{i=1}^{n-1} \Big(
 \sum_{r=1}^{p-1} \frac{H_{i-1}^{(p-r+1)}}{i^r} +2 \frac{H_{i-1}}{i^p} \Big).
$$
\end{lem}
\begin{proof}
  lhs $= \sum_{i=1}^{n-1} \frac{H_{n-i-1}}{i^p} = \sum_{i=1}^{n-1}
  \sum_{j=1}^{n-i-1} \frac{1}{i^p j} = \sum_{j=1}^{n-1}
  \sum_{i=1}^{j-1}  \frac{1}{i^p (j-i)} = $ rhs, using the partial
  fractions $ \frac{1}{i^p (j-i)} = \sum_{r=1}^{p} \frac{1}{i^r
  j^{p-r+1}} +\frac{1}{(j-i) j^p}$.
\end{proof}

\begin{lem} With $ b_i := \sum_{r=1}^i \frac{(-1)^r B_r}{r}$, we have for
  $|x|<1$:
  \begin{align}
    \mathfrak{S}_{1,1}(x) &= \Li_2(x) = \sum_{i\geq 1} \bigg( \frac{(-1)^i
  B_i}{i^2} x^i -\frac{b_i}{i} \Big(\frac x{x-1}\Big)^i \bigg)\label{eq:S11}\\
    \mathfrak{S}_{1,2}(x) &= \sum_{n\geq 1} \frac{H_{n-1}}{n^2} x^n
    = \sum_{i\geq 1} \bigg( \frac{(-1)^i B_i H_{i-1}}{i^2} x^i
-\frac{b_i}{i} [H_{i-1}+\Li_1(x)] \Big(\frac x{x-1}\Big)^i \bigg)\label{eq:S12}
  \end{align}
\end{lem}
\begin{proof}
  The first equation is equivalent to an identity for Bernoulli
  numbers:
$$
\sum_{r=1}^{n-1} \frac{(-1)^r B_r}{r} \sum_{l=r}^n \frac{(-1)^l}{l}
\binom{n-1}{n-l} = -\frac{1}{n^2}  \qquad (n\geq 2)
$$
in which the second sum equals $\frac{1}{n} \sum_{l=r}^n (-1)^l
\binom{n}{l} = -\frac{1}{n} \sum_{l=0}^{r-1} (-1)^l \binom{n}{l}
=-\frac{(-1)^r}{n} \binom{n-1}{r-1} =-\frac{(-1)^r r}{n^2}
\binom{n}{r}$.  Thus the equation boils down to the well-known 
identity for Bernoulli numbers: $ \sum_{r=0}^n \binom{n}{r} B_r =B_n$
or mnemotechnically: $(B+1)^n = B_n$ (upon replacing $B^k$ by $B_k$).

The second equation is equivalent to 
$$
-\sum_{r=1}^{p-1} \frac{(-1)^r B_r}{r} \sum_{l=r}^p \frac{(-1)^l}{l}
\binom{p-1}{p-l} H_{l-1} -\sum_{n=1}^{p-1}\bigg( \sum_{r=1}^n
\frac{(-1)^r B_r}{r} \sum_{l=r}^n \frac{(-1)^l}{l} \binom{n-1}{n-l} \bigg)
\frac{1}{p-n} =\frac{H_{p-1}}{p^2}
$$
or
$$
-\sum_{r=1}^{p-1} \frac{(-1)^r B_r}{r} \sum_{l=r}^p \frac{(-1)^l}{l}
\binom{p-1}{p-l} H_{l-1}  -\sum_{n=1}^{p-1}
\frac{(-1)^n B_n}{n^2 (p-n)} +\frac{H_{p-1}^{(2)}}{p}
+2\frac{H_{p-1}}{p^2} =\frac{H_{p-1}}{p^2}
$$

Note also that 
$$ 
\sum_{l=1}^{p+1} (-1)^l \binom{p+1}{l} H_{l-1} =-\sum_{l=0}^{p}
\frac{(-1)^{l}}{l} \binom{p}{l} = H_{p}
$$
where the second equality follows by induction, while the first comes from
$$
{\rm lhs} =(-1)^{p-1} \sum_{l=0}^{p}
(-1)^{p-l} \tbinom{p+1}{p-l} H_{l} = (-1)^{p-1} \Big[ (1-x)^{p+1}
\frac{-\log(1-x)}{1-x} \Big]_{x^p} = (-1)^{p-1} \sum_{l=0}^{p}
\frac{(-1)^{p-l}}{l} \binom{p}{p-l}.
$$
Hence the second equality is equivalent to the following conjecture,
which has resisted the author's best efforts.
\end{proof}
\begin{conj} \label{conj:bernoulli} For $p$ a positive integer:
 $$
\sum_{r=1}^{p-1} \frac{(-1)^r B_r}{r} \Big( \sum_{l=r}^p (-1)^l
\binom{p}{l} H_{l-1} +\frac{1}{r}+\frac{1}{p-r} \Big) = H_{p-1}^{(2)}
+\frac{1}{p} H_{p-1}.
$$
\end{conj}

\section*{Conclusion}

Though an application of the Euler--Maclaurin formula is nothing
distinguished, it turns out that its use for the four sums of
sections~4,5,6,7 resp. brings along a wealth of by-products about
Stirling numbers, their relation to harmonic numbers, their
asymptotics, about mathematical constants and their representation as formal
(diverging) sums over rational (Bernoulli) numbers, about more general
asymptotics of complex functions (incl. incomplete gamma functions),
as well as algebraic manipulations on polylogs and Nielsen integrals.
All this research was only possible because we used the asymp$_k$
trick and numerical mathematics. 

Extensions of this paper would be doing the same for other sums
containing logarithms, foremost $\sum k^q/(\log k)^p$, ~$\sum (n-k)^q
(\log k)^p$, but we do not expect any new property.  Perhaps only more
unknown mathematical constants would come to light.  Otherwise, the
constants $C_{p,q}$ of section~7 still await an exact form.  Further,
one could attempt to prove conjecture~\ref{conj:bernoulli} or write
down the similar (and more complex) identities that one obtains when
carefully going through the proof of lemma~\ref{logkp-(n-k)}.  Perhaps
even more bizarre identities would show up by carefully analyzing
what happens in the proofs of lemmas~\ref{logk-(n-k)q} and
\ref{logkp-(n-k)q}.

\subsection*{Acknowledgments}

It is a pleasure to acknowledge fruitful discussions with Johannes
Bl\"umlein, Karl Dilcher, Steve Finch, Philippe Flajolet, Herbert
Gangl, Hsien-Kuei Hwang, Pieter Moree, Boris Moroz, Robert Osburn,
Yiannis Petridis, Simon Plouffe, Zhi-Wei Sun, Nico Temme and Don
Zagier.  This project was supported by the Max-Planck-Institute in
Bonn.

\appendix
\section{The asymp$_k$ trick}

Assume we are given numerically the first hundreds of terms of a
converging sequence $(s_n)$, $(n\in\N)$, and that its asymptotic
expansion goes in inverse powers of $n$, ie. $s:= c_0 + {c_1 \over
  n}  +{c_2 \over n^2}  +\dots  $. \\
Goal: determine the coefficient $c_0$.\\
Trick (by Don Zagier): apply the operator ${1\over k!} \p^k n^k $ on
$s$ ($k\in\N$) to find
$$
c_0 +(-1)^k {c_{k+1} \over n^{k+1}} + \dots + (-1)^{k+l} {k+l-1
  \choose l-1} {c_{k+l} \over n^{k+l}} + \dots
$$
Hence this gives $k$ more digits of precision for $c_0$, as long as $k$
is not too big (ie. the binomials not too big).  Call this operation
\underline{asymp$_k$}.  In practice, the operator $\p$ is the
difference operator $\Delta s:= s_{n+1} -s_n$.

To determine $c_1$, subtract $c_0$, multiply by $n$ and apply
asymp$_k$ or: differentiate (ie. take successive differences) and
multiply by $-n^2$. 

The crucial point in the success of this trick is that the errors
generated by the difference operator on a monomial are
themselves monomials of lower powers:  $\Delta \frac{1}{n^j} =
\frac{1}{(n+1)^j} - \frac{1}{n^j} = \frac{-j}{n^{j+1}}
+\frac{j(j+1)/2}{n^{j+2}} +\cdots$.  These will be swept away at the next
applications of $\Delta$.  The same would not be true if the operator
acted on terms like $(\log n)^j$.

NB: we can see whether the $k$ decimals are correct by checking the
convergence of the series: if we have 400 terms, say, write every
80th term in a column (ie. 5 terms in total) and see how quickly the
digits agree.


\end{document}